\documentclass[12pt]{article}
\author{Yuval Peres\footnote{Research supported
in part by NSF grants
\#DMS-0104073 and \#DMS-0244479.}  \ and B\'alint
Vir\'ag\footnote{Research supported in part by NSF
grant \#DMS-0206781.}}
\title{Zeros of the i.i.d.\ Gaussian power series: \\
a conformally invariant determinantal process}
\date{January 20, 2005}

\usepackage{amsfonts}
\usepackage{graphicx}
\usepackage{amsmath}
\usepackage{amsthm}

\newtheorem{theorem}{Theorem}

\newtheorem{proposition}[theorem]{Proposition}
\newtheorem{lemma}[theorem]{Lemma}
\newtheorem{corollary}[theorem]{Corollary}

\newtheorem{fact}[theorem]{Fact}
\newtheorem{protoeg}[theorem]{Example}
\newtheorem{protoremark}[theorem]{Remark}
\newtheorem{protodefinition}[theorem]{Definition}

\newenvironment{remark}{\begin{protoremark}\rm}{\end{protoremark}}

\renewenvironment{proof}{\par \trivlist
 \itemindent\parindent \item[\hskip\labelsep\sc Proof.]
 \ignorespaces}{\qed\endtrivlist}
\newenvironment{proofof}[1]{\par \trivlist
 \itemindent\parindent \item[\hskip\labelsep\sc Proof of #1.]
 \ignorespaces}{\qed\endtrivlist}

\newcommand\mnote[1]{} 
\newcommand\be{\begin{equation}}
\newcommand\bel[1]{{\mnote{#1}}\be\label{#1}}
\newcommand\ee{\end{equation}}
\newcommand\lb[1]{\label{#1}\mnote{#1}}
\newcommand{\PP}{\mathcal P}
\newcommand{\NN}{\mathcal N}
\newcommand{\comment}[1]{}
\newcommand{\eps}{\varepsilon}

\newcommand{\UU}{{\mathbb U}}
\newcommand{\R}{{\mathbb R}}
\newcommand{\CC}{{\mathbb C}}

\newcommand{\ev}{\mbox{\bf E}}
\newcommand{\pr}{\mbox{\bf P}}
\newcommand{\one}{{\mathbf 1}}
\newcommand{\as}{\mbox{\hspace{.3cm} a.s.}}
\newcommand{\dist}{\mbox{\rm dist}}

\newcommand{\perm}{\mbox{\rm perm}}
\newcommand{\Poi}{\mbox{\rm Poi}}
\newcommand{\Var}{\mbox{\rm Var}}
\newcommand{\Vol}{\mbox{\rm Vol}}

\newcommand{\Cov}{\mbox{\rm Cov}}
\newcommand{\sm}{{\raise0.3ex\hbox{$\scriptstyle \setminus$}}}
\newcommand{\intc}{\int_0^{2\pi}}

\newcommand{\re}[1]{(\ref{#1})}
\newcommand{\gi}{\,|\,}
\newcommand{\area}{\operatorname{area}}
\newcommand{\szego}{Szeg\H o}
\newcommand{\h}{h}

\newcommand{\wf}{\widetilde{f}}
\newcommand{\wN}{\widetilde{N}}
\newcommand{\G}{G}
\newcommand{\bb}{\beta}
\newcommand{\hL}{\widehat{L}}
\newcommand{\hNN}{\widehat{\NN}}
\newcommand{\zb}{\overline{z}}
\newcommand{\ed}{\stackrel{d}{=}}

\def\summ{{\sum\limits}}
\def\intt{{\int\limits}}
\def\prodd{{\prod\limits}}
\def\mb{\mbox}
\def\l{\left}
\def\r{\right}
\def\sig{\sigma}

\def\eps{\epsilon}

\usepackage{natbib}

\begin{document}
\maketitle
\begin{abstract}
 Consider the zero set of the random power series
$f(z)=\sum a_n z^n$ with i.i.d.\ complex Gaussian
coefficients $a_n$. We show that these zeros form a
determinantal process: more precisely, their joint
intensity can be written as a minor of the Bergman
kernel.  We show that the number of zeros of $f$ in a
disk of radius $r$ about the origin has the same
distribution as the sum of independent
$\{0,1\}$-valued random variables $X_k$, where
$P(X_k=1)=r^{2k}$. Moreover, the set of absolute
values of the zeros of $f$ has the same distribution
as the set $\{U_k^{1/(2k)}\}$ where the $U_k$ are
i.i.d.\ random variables uniform in $[0,1]$. The
repulsion between zeros can be studied via a dynamic
version where the coefficients perform Brownian
motion; we show that this dynamics is conformally
invariant.
\end{abstract}

\section{Introduction}
Consider the random power series \bel{series} f_{\UU}(z)
=\sum_{n=0}^\infty a_n z^n \ee where $\{a_n\}$ are independent
standard complex Gaussian random variables (with density $
e^{-z\overline z}  /\pi.$) The radius of convergence of the series
is a.s.\ 1, and the set of zeros forms a point process $Z_\UU$ in
the unit disk $\UU$.
 Zeros of Gaussian power series have been studied starting with
 \cite{offord65},
  since these series are limits of random Gaussian polynomials.
  In the last decade,
physicists have introduced  a new perspective, by interpreting the zeros
of  a Gaussian polynomial  as a  gas of interacting particles,
see \cite{hannay96}, \cite{leboeuf99} and the references therein.
Much of the recent interest in Gaussian analytic
functions was spurred by the papers
\cite{ek95} and \citet*{bsz00}.
A fundamental
property of $Z_\UU$ is the invariance of its distribution
under M\"obius transformations that preserve the unit
disk; see Section 2 for an explanation,  and \cite{ts03} for
references.

Our main new discovery is that the zeros $Z_{\UU}$ form a {\em
determinantal process\/},  and this yields an explicit formula
for the distribution of the number of zeros in a disk.
Furthermore, we show that the process $Z_{\UU}$ admits a
conformally invariant evolution which elucidates the repulsion
between zeros.

Given a random function $f$ and points $z_1,\ldots, z_n$, let
$p_\eps(z_1,\ldots, z_n)$ denote the probability that for all $1
\le i \le n$, there is a zero of $f$ in the disk of radius $\eps$
centered at $z_i$. The {\bf joint intensity} of the point process
of zeros of $f$, also known as the $n$-{\em point correlation
function}, is defined by the limit \bel{jointdef}
p(z_1,\ldots,z_n) =\lim_{\eps\to 0}
  \frac{p_\eps(z_1,\ldots,z_n)}{\pi^n\eps^{2n}}
\ee
when it exists; see (\ref{soshcor1}) for a related integral formula.

\begin{theorem}\label{detformula}
 The joint intensity of zeros for the
 i.i.d.\ Gaussian power series (\ref{series}) in the
 unit disk
 exists, and satisfies
\bel{det1} p(z_1,\ldots,z_n)=\pi^{-n}\det
\left[\frac{1}{(1-z_i\overline z_j)^{2}}\right]_{i,j}.
\ee
\end{theorem}
Thus the zero set of the i.i.d.\ series $f_\UU(z)$ is a
determinantal process in $\UU$, governed by the {\bf Bergman
kernel} $K_\UU(z,w)=\pi^{-1}(1-z\overline w)^{-2}$; see
\cite{soshnikov00} for a survey of determinantal processes. In
particular, (\ref{det1}) extends the known fact that
$p(z_1,z_2)<p(z_1)p(z_2)$ for all $z_1,z_2 \in \UU$, i.e., the
zeros are negatively correlated.
In fact,
 $Z_\UU$ is the only process of zeros of a
Gaussian analytic function which is negatively correlated and has
a M\"obius invariant law; see Section \ref{s.ci}.

The determinant formula for the joint intensity allows us to
determine the distribution of the number of zeros of $f_\UU$ in a disk,
and identify the law of the moduli of the zeros.
\begin{theorem}\lb{distribution}
\begin{description}
\item{\bf (i)} The number  $N_r=|Z_\UU \cap B_r(0)|$ of zeros of $f_\UU$ in the
disk of Euclidean radius $r$ about $0$, satisfies
\bel{laplace}
\ev (1+s)^{N_r}=\prod_{k=1}^\infty
(1+r^{2k}s)
\ee
for all real $s$. Thus $N_r$ has the
same distribution as $\sum_{k=1}^\infty X_k$ where
$\{X_k\}$ is a sequence of independent
$\{0,1\}$-valued random variables with
$\pr(X_k=1)=r^{2k}$.
\item{\bf (ii)} Moreover, the set of moduli
 $\{|z| : f_\UU(z)=0\}$
 has the same law as the set $\{U_k^{1/(2k)}\}$, where $\{U_k\}$ are
i.i.d.\ random variables uniform in $[0,1]$.
\end{description}
\end{theorem}

From Theorem \ref{distribution} we readily obtain the
asymptotics of the hole probability $\pr(N_r=0)$.
Furthermore, the infinite product in (\ref{laplace})
occurs in one of Euler's partition identities, see
(\ref{euler}), and this connection yields part (ii) of
the next corollary.
\begin{corollary}\lb{hole}
\begin{description}
\item{\bf (i)}
Let $h=4\pi r^2/(1-r^2)$, the hyperbolic area of $B_r(0)$. As $r
\uparrow 1$,  we have
$$
\pr(N_r=0)=\exp\Bigl({\frac{-\pi h +o(h)}{24}}\Bigr)=
\exp \Bigl({\frac{-\pi^2+o(1)}{12(1-r)}} \Bigr)\, .
$$
\item{\bf (ii)} The binomial moments of $N_r$ equal
 $$
 \ev \binom{N_r}{k}
=\frac{r^{k(k+1)}}{(1-r^2)(1-r^4)\cdots (1-r^{2k})}.
 $$
\item{\bf (iii)} The ratio
$(N_r-\mu_r)/\sigma_r$ converges in law to
standard normal as $r \uparrow 1$, where
$$\mu_r=\ev N_r = \frac{r^2}{1-r^2} ,\quad \mbox{ \rm and } \sigma^2_r=
\Var\, N_r=\frac{r^2}{1-r^4} \, .$$
\end{description}
\end{corollary}

\subsection{General domains} The covariance structure $\ev
\Bigl(f_\UU(z)\overline{
f_\UU(w)}\Bigr)=(1-z\overline w)^{-1}$ 
equals  $2\pi$ times the {\bf \szego\ kernel}
$S_\UU(z,w)=(2\pi)^{-1}(1-z\overline w)^{-1}$ in the
unit disk. The \szego\ kernel $S_D(z,w)$ and the
Bergman kernel $K_D(z,w)$ are defined, and positive
definite, for any bounded  planar domain $D$ with a
smooth
boundary. (See the next section or \cite{bell} 
for information on the \szego\ and Bergman kernels.) For such
domains we can consider the Gaussian analytic function $f_D(z)$
with covariance structure $2 \pi S_D$ in $D$ (an explicit formula
for $f_D$ is given in (\ref{genf})). Recall
 that a {\bf Gaussian analytic function} in $D$ is a random analytic
  function $f$ such that for any choice of $z_1,\ldots,z_n$ in $D$,
the random variables $f(z_1),\ldots,f(z_n)$ have complex Gaussian
joint distribution.

\begin{corollary}\label{cdetformula}
 Let $D$ be a simply connected  bounded planar domain,
with a $C^\infty$ smooth boundary.
 The joint intensity of zeros for the
Gaussian analytic function $f_D$ is given by the
determinant of the Bergman kernel
$$
p(z_1,\ldots,z_n)=\det [K_D(z_i,z_j)]_{i,j}.
$$
\end{corollary}
Note that for simply connected domains as in the
corollary, the Bergman and \szego\ kernels satisfy
$K_D(z,w)= 4\pi S_D(z,w)^2$, see \cite{bell}, Theorem
23.1.

\subsection{The one-parameter family of M\"obius-invariant zero sets}
For $\rho>0$, let $Z_{\UU,\rho}$ denote the zero set of
\bel{genrho} f_{\UU,\rho}(z)
 =  \sum_{n=0}^\infty
 \binom{-\rho}{n}^\frac{1}{2} a_n
 z^n \, ,
\ee
where $\{a_n\}$ are i.i.d.\ standard complex Gaussians.
In particular, $f_{\UU,1}$ has the same distribution as $f_\UU$.
As explained in \cite{ts03} (see also \cite{br02}), for any $\rho>0$,
the distribution of  $Z_{\UU,\rho}$
 is invariant under M\"obius transformations that
 preserve $\UU$. Moreover, these are the only zero sets of Gaussian
analytic functions with this invariance property.
However, only $\rho=1$ yields a determinantal zero process.

\medskip
Taking $n=1$ in Theorem \ref{detformula}, one recovers the
well-known formula  $(1-|z|^2)^{-2}/\pi$ for the intensity of
$Z_\UU$. More generally,
the intensity of $Z_{\UU,\rho}$ in
$\UU$ is $\rho/[\pi(1-|z|^2)^2]$, see \cite{sodin00}.
  It follows that the
expected number of zeros in a Borel set  $\Lambda \subset \UU$ is
$\rho/(4\pi)$ times the hyperbolic area $A(\Lambda)=\int_\Lambda
\frac{4dz}{(1-|z|^2)^2}$. (Integration is with respect to planar
Lebesgue measure.) This can also be inferred from Proposition
\ref{permformula} below. In Section \ref{s.lln}, we prove the following
 law of large numbers. 
\begin{proposition}\label{lln}
Let $\rho>0$, and suppose that $\{\Lambda_h\}_{h>0}$ is
an increasing family of Borel sets
in $\UU$, parameterized by hyperbolic area $h=A(\Lambda_h)$. Then
the number  $N(\h)=|Z_{\UU} \cap \Lambda_h|$ of zeros of $f_{\UU,\rho}$
in $\Lambda_h$ satisfies
$$
\lim_{h\to \infty} \frac {N(\h)}{h} =\frac{\rho}{4\pi} \as
$$
\end{proposition}

\subsection{Reconstruction of $|f_{\UU,\rho}|$ from its zeros}
\begin{theorem} \lb{recprod}
\begin{description}
\item{\bf (i)}
Let $\rho>0$. Consider the random
function $f_{\UU,\rho}\, $, and order its zero set
$Z_{\UU,\rho}$ in increasing absolute value,
  as $\{z_k\}_{k=1}^\infty$. Then
\bel{e.recp2} |f_{\UU,\rho}(0)|=c_\rho
\prod_{k=1}^\infty e^{\rho/(2k)} |z_k| \quad a.s. \ee
where
$c_\rho=e^{(\rho-\gamma-\gamma\rho)/2}\rho^{-\rho/2}$
and $\gamma=\lim_n \Bigl(\sum_{k=1}^n \frac{1}{k}-\log
n \Bigr)$ is Euler's constant.

\item{\bf(ii)} More generally, given $\zeta \in \UU$, let
$\{\zeta_k\}_{k=1}^\infty$ be $Z_{\UU,\rho}$, ordered
in increasing hyperbolic distance from $\zeta$. Then
\bel{e.recp3}
|f_{\UU,\rho}(\zeta)|=c_\rho(1-|\zeta|^2)^{-\rho/2}
\prod_{k=1}^\infty e^{\rho/(2k)}
\Bigl|\frac{\zeta_k-\zeta} {1-\overline{\zeta}
\zeta_k} \Bigr| \,. \ee
\end{description}
Thus  the analytic function
$f_{\UU,\rho}(z)$ is determined by its zero set, up to
multiplication by a constant
of modulus 1.
\end{theorem}
This theorem is proved in Section \ref{s.rec}.

\subsection{Dynamics}
In order to understand the negative correlations for zeros of
$f_\UU$, we consider a dynamic version of the zero set $Z_\UU$.
Denote by $Z_{\UU}(t)$ the zero set of the power series
$\sum_{n=0}^\infty a_n(t)z^n$, where the coefficients $a_n(t)$ are
independent stationary complex Ornstein-Uhlenbeck processes; in
other words, $a_n(t)=e^{-t/2} W_n(e^t)$, where $\{W_n(\cdot)\}_{n
\ge 0}$ are independent complex Brownian motions.

A direct calculation gives that for the process
$Z_{\UU}$, the intensity ratio
$p(z_1,z_2)/[p(z_1)p(z_2)]$ is strictly less then 1
and decreases to 0 with as the hyperbolic distance
between $z_1$ and $z_2$ tends to 0. This repulsion
suggests
 that when two zeros get close,
there is a drift in their motion that pushes them apart. However,
this is not the case. Instead, we have the following.
\begin{theorem} \lb{thm:dyn}
Consider the process of zeros $\{Z_{\UU}(t)\}$ in the
unit disk and condition on the event that  at time $t=0$ there is a
zero at the origin, i.e., $0\in Z_\UU(0)$. The movement of this zero is
then described by an SDE which at time $t=0$ has the
form
$$
dz= \sigma dW,
$$
where $W$ is complex Brownian motion, there is no
drift term, and
$$
1/\sigma = |f'_\UU(0)|=   c_1 \prod_{k=2}^\infty
e^{1/(2k)} |z_k| \quad a.s.
$$
Heuristically, any zero of $f_\UU$  oscillates faster when there are other
zeros nearby; this causes repulsion.
\end{theorem}

Analogous processes $Z_D(t)$ can be defined in general
domains, and we shall show in Section \ref{s.dyn} that
the family of processes $Z_D(t)$ is conformally
invariant (no time change is needed). Theorem \ref{thm:dyn}
can be extended to $\rho \ne 1$ as well.

\medskip \noindent {\bf Conditioning to have a zero at a
 given location.} It
is important to note that the distribution of $f_\UU$
given that its value is zero at 0 is {\it different}
from the conditional distribution of $f_\UU$ given that its zero
set has a point at 0. In particular, in the second
case the conditional distribution of the coefficient
$a_1$ is not Gaussian. The reason for this is that the two ways of
conditioning are defined by the limits as $\eps\to 0$ of
two different conditional distributions. In the first
case, we condition on $|f_\UU(0)|<\eps$. In the second, we
condition on $f_\UU$ having a zero in the disk $B_\eps(0)$
of radius $\eps$ about 0; the latter conditioning affects
the distribution of $a_1$. See Lemma \ref{condz} in Section \ref{s.rec}.

\subsection{Hammersley's formula}
The starting point of the proof of Theorem
\ref{detformula} is a general permanent formula for
the joint intensity of zeros for Gaussian analytic
functions. A version for polynomials is due to
\cite{hammersley56} and \cite{friedman90}. The
permanent form \re{density} for Gaussian polynomials
appears in the physics literature (\cite{hannay96}).
Closely related formulae for correlations between
zeros of random sections of a positive holomorphic
line bundle over a compact complex manifold were
established by \citet*{bsz00}.

The version we need is for Gaussian analytic
functions, that are not necessarily polynomials.

\begin{proposition}\lb{permformula}
Let $f$ be a Gaussian analytic function
 in a planar domain $D$ such that $\ev f(z)=0$ for all $z \in D$.
Given points
$z_1,\ldots, z_n\in D$,  consider the matrices
$$ A=\Big(\ev
f(z_i)\overline{f( z_j)}\Big); \quad
 B=\Big(\ev f'(z_i)\overline{f(z_j)}\Big); \quad
 C=\Big(\ev f'(z_i)\overline{f'(z_j)}\Big).
$$
Assume that $A$ is nonsingular.

\noindent{\bf (i)}  The  joint intensity for the
zeros of $f$ exists and satisfies
%
\bel{prehan}
p(z_1,\ldots, z_n) = \frac{\ev\Bigl(|f'(z_1)\cdots
f'(z_n)|^2 \, \Big| \, f(z_1)=\cdots =f(z_n)=0\Bigr)}{\det (\pi
A)}.
\ee
Consequently,
\bel{density}
p(z_1,\ldots,z_n)=\frac{\perm(C-BA^{-1}B^*)}{\det(\pi
A)}. \ee
 \noindent{\bf (ii)}
Assume that $A=A(z_1,\ldots, z_n)$ is nonsingular when
$z_1,\ldots,z_n\in D$ are distinct. Let $Z_f^{\wedge
n}$ denote the set of  $n$-tuples of distinct zeros of
$f$. Then for any Borel set $B\subset D^n$ we have
 \bel{soshcor1} \ev\,
\#(B\cap Z^{\wedge n}) =\int_Bp(z_1,\dots
,z_n)\,dz_1\dots dz_n \,. \ee
\end{proposition}
The proof of this proposition is given in Section
\ref{general}.

In the derivation of Theorem \ref{detformula} from
proposition \ref{permformula}, we use conformal
invariance, the i.i.d.\ property of the coefficients,
and the beautiful determinant-permanent identity
(\ref{borchardt}) of \cite{borchardt1855}.

\medskip
\noindent{\bf Remarks on the literature}. A nice
introduction to the theory of Gaussian analytic
functions is given in \cite{sodin00}; for earlier
results, see  \cite{hammersley56}, 
\cite{friedman90}, \cite{bbl92}, \cite{kostlan93},
\cite{ek95} and \cite{hannay96}. Close to the topic of
this paper are \cite{sz03},
\cite{ts03}.
Determinantal processes are also being intensively
studied,
see \cite{soshnikov00}. Theorem \ref{detformula} provides
further evidence for the analogy, suggested by
\cite{leboeuf99}, between zeros of Gaussian
polynomials and the Ginibre ensemble of eigenvalues of
(non-hermitian) random matrices with i.i.d.\ Gaussian
entries, which is known to be determinantal.

\section{Conformal invariance and preliminaries}\lb{s.ci}
\noindent {\bf Complex Gaussian random variables.}
Recall that a standard complex Gaussian random
variable $a$ has density $ e^{-z\overline z} /\pi$,
expected value $0$, and variance $\ev a
\overline{a}$=1. A vector $V$ of random variables has a
complex Gaussian (joint) distribution if
there is a determinstic vector $V_0$ such that $V-V_0$ is
the image under a linear map
of a vector of i.i.d.\ standard complex Gaussian random
variables.

If $X$, $Y$ are real Gaussian random variables of mean zero, then
$X+iY$ is complex Gaussian if and only if $X$, $Y$ are
independent and have the same variance. A complex
Gaussian random variable $Z$ with $\ev Z=0$ satisfies
 $\ev Z^n=0$ for any integer $n \ge 1$.

\medskip
\noindent {\bf The complex Gaussian power series.}
Recall the power series $f_{\UU}$ in (\ref{series}). A
 Borel-Cantelli argument shows that the radius
of convergence of $f_{\UU}$ equals 1 a.s. Clearly, the
joint distributions of $f_{\UU}(z_k)$ for any finite
collection $\{z_k\}$ are complex Gaussian, so the
values of $f_{\UU}$ form a complex Gaussian ensemble.
Since $f_{\UU}$ is continuous, its distribution is
determined by the covariance structure
 \bel{covariance}
\ev \Bigl(f_{\UU}(z) \overline {f_{\UU}(w)} \Bigr) =
  \sum_{n=0}^\infty (z \overline{w})^n = (1 - z\overline{w})^{-1}.
 \ee
The right hand side is $2 \pi$
times the \szego\ kernel in the unit disk; it suggests a natural
way to generalize the power series $f_\UU$.

\medskip
\noindent {\bf The \szego\ kernel.} Let $D$ be a
bounded planar domain with a $C^\infty$ smooth boundary
(the regularity
assumption can be weakened).
%
Consider the set of complex analytic functions in $D$
which extend continuously to the boundary $\partial
D$. The classical Hardy space
$H^2(D)$ is given by the $L^2$-closure of this set
with respect to length measure on $\partial D$.
Every element of $H^2(D)$ can be identified with a
unique analytic function in $D$ via the Cauchy
integral (see \cite{bell}, Section 6).

Consider an orthonormal basis $\{\psi_n\}_{n\ge 0}$ for $H^2(D)$;
e.g. in the unit disk, take $\psi_n(z)=\frac{z^n}{\sqrt{2\pi}}$
for $n \ge 0$.
  Use i.i.d.\ complex Gaussians $\{a_n\}_{n\ge 0}$
 to define the random analytic
  function
 \bel{genf}
f_D(z)=\sqrt{2\pi}\sum_{n=0}^{\infty} a_n \psi_n(z) \,.
 \ee
(cf.\ (6) in \cite{sz03}).
The factor of $\sqrt{2\pi}$ is included just to simplify formulas
in the case where $D$ is the unit disk. The covariance function of
$f_D$ is given by $2 \pi S_D(z,w)$, where
 \bel{szdef} S_D(z,w)=\sum_{n=0}^{\infty}
\psi_n(z)\overline{\psi_n(w)}
 \ee
 is the {\szego\ kernel} in $D$.
The \szego\ kernel $S_D$ does not depend on the choice of
 orthonormal basis and is positive
definite (i.e. for points $z_j\in D$ the matrix
$(S_D(z_j,z_k))_{j,k}$ is positive definite).

Let $T:\Lambda\to D$ be a conformal homeomorphism between
two bounded domains with $C^\infty$ smooth boundary. The
derivative $T'$ of the conformal map has a
well-defined square root, see
\cite{bell} p.\ 43. If $\{\psi_n\}_{n\ge 0}$ is an
orthonormal basis for $H^2(D)$, then
$\{ \sqrt{T'} \cdot (\psi_n\circ  T)  \}_{n\ge 0}$ forms an orthonormal
basis for $H^2(\Lambda)$. In particular, the \szego\ kernels
satisfy
 \bel{szegoci} S_\Lambda(z,w)=T'(z)^{1/2}\overline{
T'(w)^{1/2}}S_{D}(T(z),T(w)).
 \ee
When $D$ is a simply connected domain, it follows from
the transformation formula \re{szegoci} that $S_D$
does not vanish in the interior of $D$, so for
arbitrary $\rho>0$ powers $S_D^\rho$ are defined.

Let $\{\eta_n\}_{n\ge 0}$ be an orthonormal basis of the
subspace of complex analytic functions in $L^2(D)$ with respect to
Lebesgue area measure. The Bergman kernel
 $$ K_D(z,w)=\sum_{n=0}^{\infty}
\eta_n(z)\overline{\eta_n(w)}
 $$
is independent of the basis chosen, see \cite{nehari}, formula
(132).

\medskip
\noindent{\bf The \szego\ random functions with parameter $\rho$.}
Recall the one-parameter family of Gaussian analytic functions
$f_{\UU,\rho}$ defined in (\ref{genrho}).
The binomial expansion yields that the covariance structure  $\ev
\Bigl(f_{\UU,\rho}(z)\overline{f_{\UU,\rho}(w)}\Bigr)$
 equals
\bel{binexp}
\sum_{n=0}^\infty \Bigl| \binom{-\rho}{n} \Bigr| \,
 z^n \overline{w}^n =\sum_{n=0}^\infty \binom{-\rho}{n}
 (-z\overline{w})^n= (1-z\overline w)^{-\rho}=[2\pi S_\UU(z,w)]^\rho \,.
\ee
The invariance of the  distribution of $Z_{\UU,\rho}$
under M\"obius transformations of the unit disk is a
special case of the following.
\begin{proposition}\lb{ci}
Let $D$ be a bounded planar domain with a $C^\infty$ boundary
and let $\rho>0$. Suppose that
either {\bf (i)} $D$ is simply connected or {\bf (ii)} $\rho$ is an integer.
Then there is a mean zero Gaussian analytic function
$f_{D,\rho}$ in $D$ with covariance structure
$$
\ev \Bigl(f_{D,\rho}(z) \overline {f_{D,\rho}(w)} \Bigr)
=[2\pi S_D(z,w)]^\rho
 \quad \mbox{ \rm for }  z,w \in D.
$$
The zero set $Z_{D,\rho}$ of $f_{D,\rho}$ has a
conformally invariant distribution:
if $\Lambda$ is another bounded domain with a smooth boundary,
 and $T:\Lambda \to D$ is a conformal
homeomorphism, then $T(Z_{\Lambda,\rho})$ has the same distribution as
$Z_{D,\rho}$. Moreover, the following two random
 functions have the same distribution:
\bel{gentran}
f_{\Lambda,\rho}(z) \stackrel{d}{=} T'(z)^{\rho/2}
\cdot(f_{D,\rho}\circ T)(z) \,.
\ee
\end{proposition}
We call the  Gaussian analytic function $f_{D,\rho}$
described in the proposition the {\bf \szego\ random
function with parameter $\rho$} in $D$.
\begin{proof}
\noindent{\bf Case (i):} $D$ is simply connected.
Let $\Psi: D \to \UU$ be a conformal map onto $\UU$,
and let $\{a_n\}$ be i.i.d.\ standard complex Gaussians.
We claim that
\bel{genrhoD} f(z)=  \Psi'(z)^{\rho/2} \sum_{n=0}^\infty
 \binom{-\rho}{n}^\frac{1}{2} a_n
 \Psi(z)^n \,
\ee
is a suitable candidate for $f_{D,\rho}$.
Indeed, repeating the calculation in (\ref{binexp}), we find that
\begin{eqnarray*}
\ev \Bigl(f(z)\overline{f(w)}\Bigr)
&=& [\Psi'(z)\overline{\Psi'(w)}]^{\rho/2}
 (1-\Psi(z)\overline{\Psi(w)})^{-\rho} \\[1ex]
&=& [\Psi'(z) \overline{\Psi'(w)}]^{\rho/2}
 \cdot[2\pi S_\UU(\Psi(z),\Psi(w))]^\rho \,.
\end{eqnarray*}
The last expression equals $[2\pi S_D(z,w)]^\rho$
by the transformation formula \re{szegoci}.
Thus we may define $f_{D,\rho}$ by the right hand side of (\ref{genrhoD}).
If $T:\Lambda \to D$ is a conformal
homeomorphism,
 then $\Psi \circ T$ is a conformal
map from $\Lambda$ to $\UU$, so (\ref{genrhoD}) and the chain rule give
the equality in law \re{gentran}. Since $T'$ does not have zeros
in $\Lambda$, multiplying $f_{D,\rho}\circ T$  by a power of $T'$
does not change its zero set in $\Lambda$, and it follows that
$T(Z_{\Lambda,\rho})$ and
$Z_{D,\rho}$ have the same distribution.
\medskip

\noindent{\bf Case (ii):} $\rho$ is an integer.
Let $\{\psi_n\}_{n\ge 0}$ be an orthonormal basis for $H^2(D)$.
  Use i.i.d.\ complex Gaussians
$\{a_{n_1,\ldots,n_\rho} \, : \, n_1,\ldots,n_\rho \ge 0 \}$
 to define the random analytic function
 \bel{genfrho2}
f_{D,\rho}(z)=(2\pi)^{\rho/2}\sum_{n_1,\ldots,n_\rho\ge 0}
 a_{n_1,\ldots,n_\rho} \psi_{n_1}(z) \cdots \psi_{n_\rho}(z)\,;
 \ee
see \cite{sodin00} for convergence. A direct calculation shows that
$f_{D,\rho}$, thus defined, satisfies
$$
\ev \Bigl(f_{D,\rho}(z) \overline{f_{D,\rho}(w)} \Bigr)
=(2\pi)^\rho\sum_{n_1,\ldots,n_\rho\ge 0}
  \psi_{n_1}(z) \overline{\psi_{n_1}(w)}
   \cdots \psi_{n_\rho}(z) \overline{\psi_{n_\rho}(w)}
=   [2\pi S_D(z,w)]^\rho \,.
$$
The transformation formula \re{szegoci} implies that
the two sides of \re{gentran} have the same covariance structure,
$[2\pi S_\Lambda(z,w)]^\rho$. This establishes \re{gentran} and completes
the proof of the proposition.
\end{proof}

The general theory of Gaussian analytic functions implies that, up
to multiplication by a deterministic analytic function, the
random functions $f_{\UU,\rho}$ are
the only Gaussian analytic functions with marginal
intensity of zeros proportional to the
hyperbolic area element.  See \cite{sodin00} for a proof.

Similarly, the zeros of the Gaussian analytic function
$$
F_{\CC,\rho}(z) = \sum_{n=0}^\infty
\left(\frac{\rho^n}{n!}\right)^\frac{1}{2} a_n z^n
$$
have distribution which is invariant under rotations and
translations of the complex plane. Note that here $\rho$ is a
simple scale parameter: $\; F_{\CC,\rho}(z) =F_{\CC,1}(\sqrt{\rho} z)$.

Letting $\rho\to\infty$ in the definition of $f_{\UU,\rho}$,  one
recovers that the limit of the rescaled point processes
$\rho^{1/2} Z_{\UU,\rho}$ is the
zero set of $F_{\CC,1}$; this phenomenon and its generalizations
have been studied by \cite{bsz00}.
\begin{figure}
\centering
\includegraphics[height=2in]{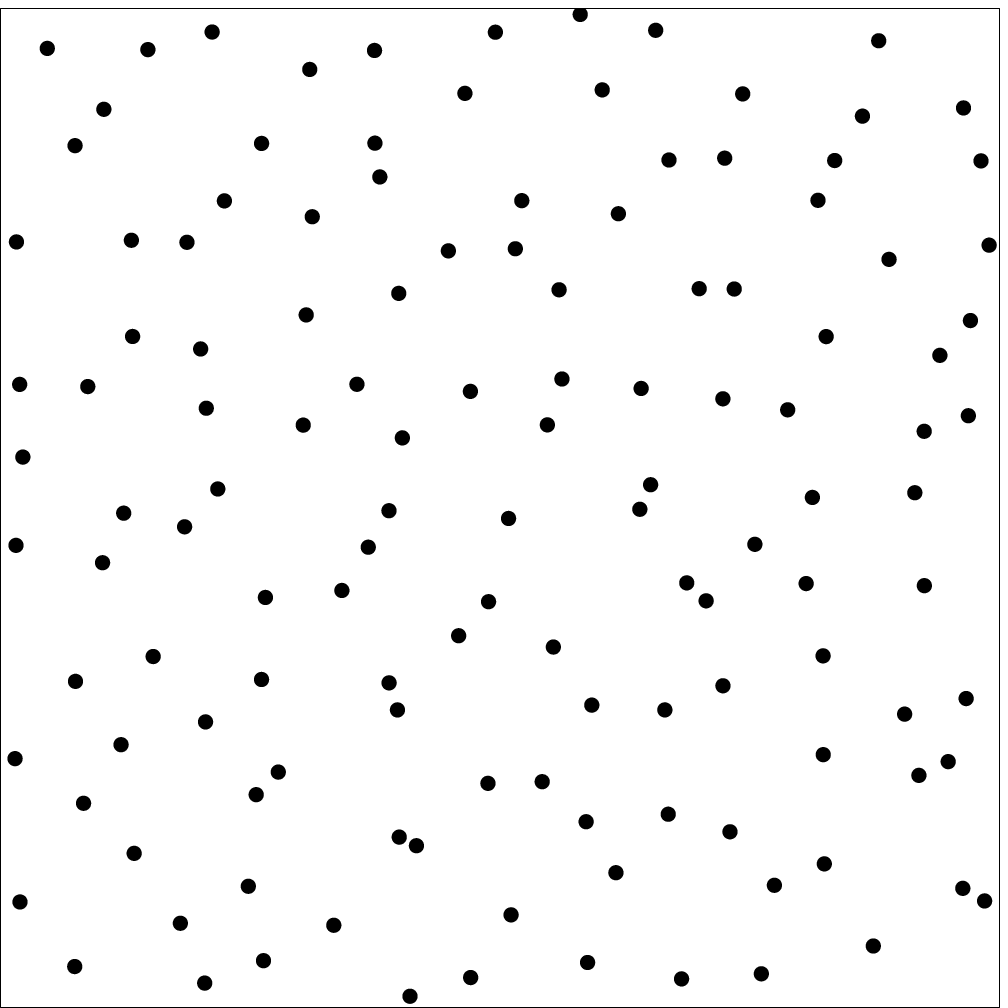}\hspace{.7in}
\includegraphics[height=2in]{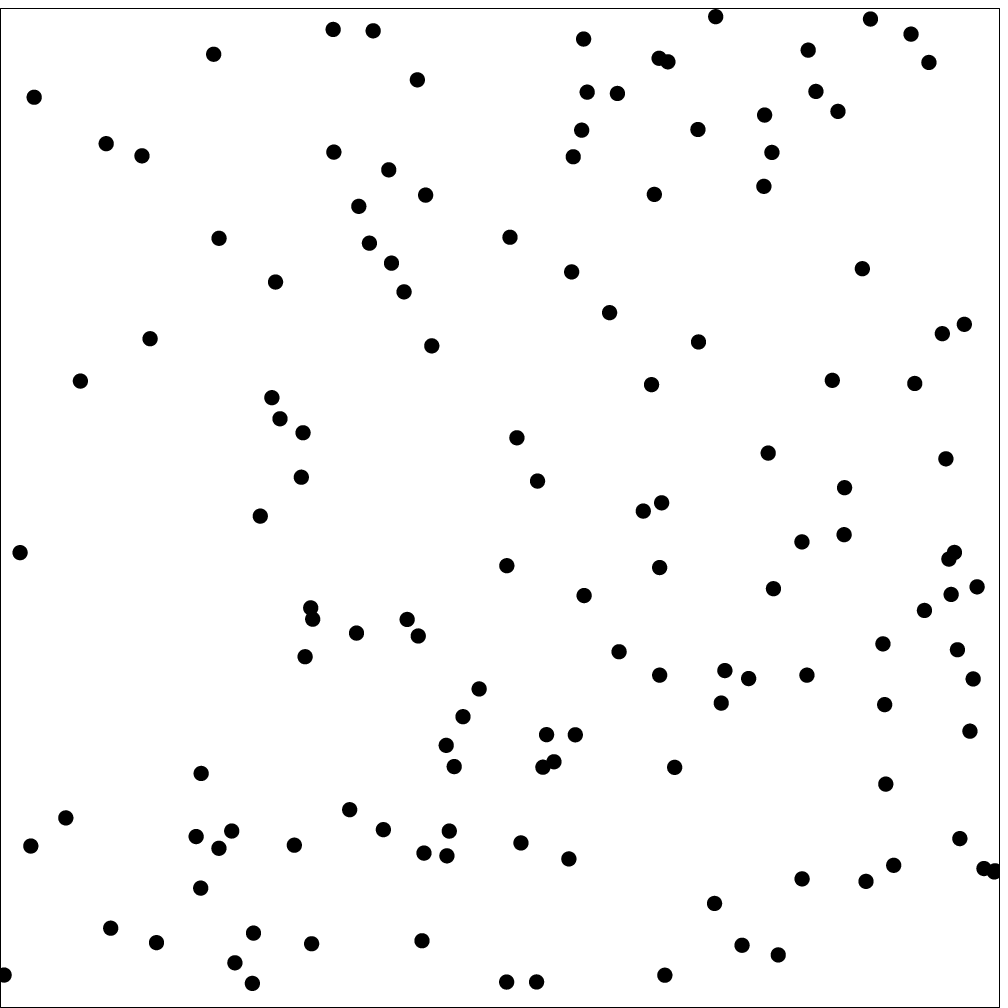}
\caption{\label{f2}The translation invariant root process and a
Poisson point process with the same intensity on the plane}
\end{figure}

\begin{figure}
\centering
\includegraphics[height=1.5in]{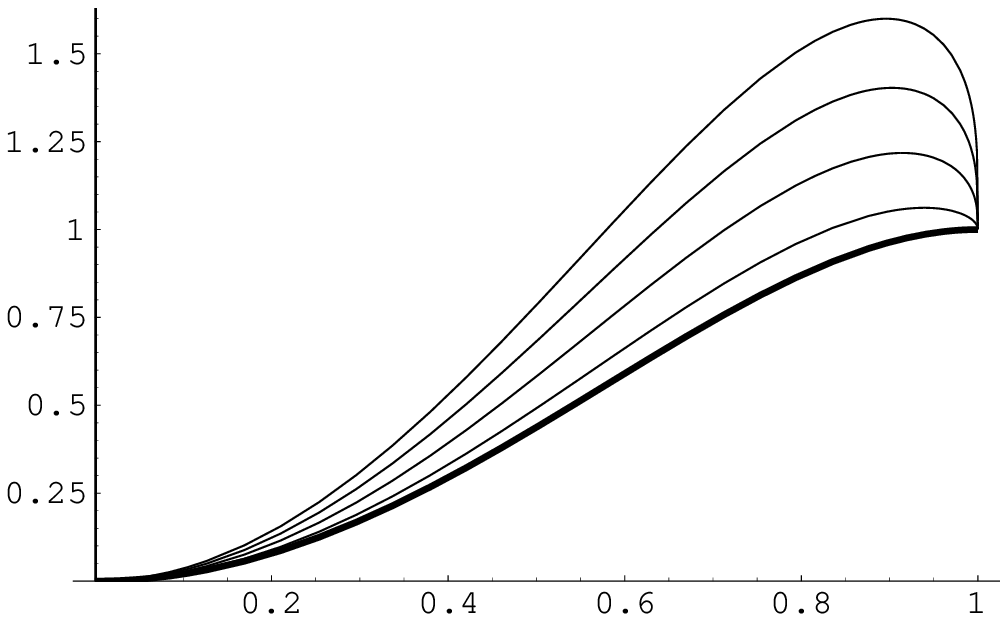}\hspace{.5in}
\includegraphics[height=1.5in]{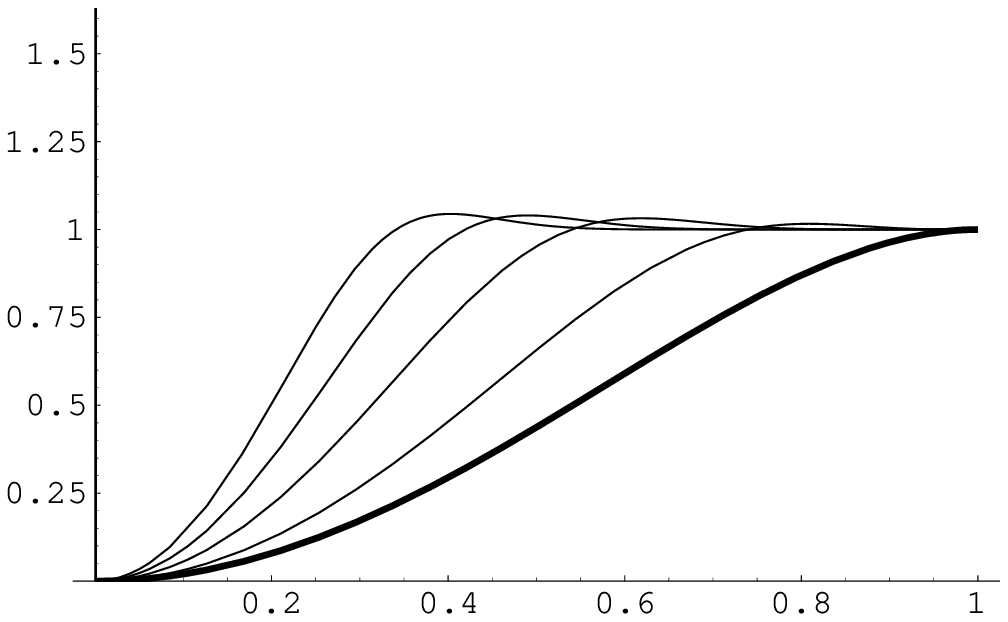}
\caption{\label{fig1}Relative intensity at $z=0$ and $z=r$ as a
function of $r$ for $\rho=1,\frac{1}{2},\frac{1}{3},\frac{1}{4
},\frac{1}{5}$ and for $\rho=1,4,9,16,25$.}
\end{figure}

Figure \ref{f2} shows a realization of the whole plane Gaussian
zero process along with a Poisson point process of the same
intensity. The orderliness of the zeros suggests that there is a
local repulsion taking place. One gets similar pictures for the
\szego\ random functions in the unit disk. The two-point intensity
for zeros at the points $r$ and $0$ is given by (\ref{density}).
The most revealing formula is the ratio $p(0,r)/(p(0)p(r))$, which
shows how far the point process is from a Poisson point process,
where this ratio is identically 1. For general $\rho$, with the
notation $s=1-r^2$ this ratio equals \bel{2point}
  \frac{1 +
  \left( {\rho ^2} - 2\,\rho  - 2 \right) \,\left( {s^\rho } + {s^{2 + 2\rho }} \right)  +
  {{\left( \rho  + 1 \right) }^2}\,\left( {s^{2\rho }} + {s^{2 + \rho }} \right)  -
  2\,{\rho ^2}\,\left( {s^{1 + \rho }} + {s^{1 + 2\rho }} \right)  + {s^{2 + 3\rho }}}
  {{{\left( 1 - {s^\rho } \right) }^3}}
 \ee
and in the case $\rho=1$ it simplifies to $
r^2(2-r^2). $ For every distance $r$, the correlation
is minimal when $\rho=1$ (see Figure \ref{fig1}). For
all values of $\rho$ different from $1$, for small
distances zeros are negatively correlated, while for
large distances the correlation is positive.

When $\rho=1$, the zeros are purely negatively correlated:  this
special phenomenon is explained by the determinantal form of the
joint intensity.

\noindent {\bf Remark.} The \szego\ random function
for $\rho=2$ $$ \sum_{n=0}^{\infty} \sqrt{n+1}a_nz^n
$$ coincides with the $n\to\infty$ limit of the logarithmic derivative of the
characteristic function of a random $n\times n$
unitary matrix (see \cite{de01}).

\medskip
 \noindent{\bf The analytic extension of white noise.} Next,
  we show that up to the
constant term, the power series $f_\UU$ has the same
distribution as the analytic extension of white noise
on the unit circle. Let $B(\cdot)$ be a standard real Brownian
motion, and let
$$
u(z) = \intc \Poi(z,e^{it}) d B(t) \,.
$$
Here the integral with respect to $B$ can be interpreted either as
a stochastic integral, or as a  Riemann-Stieltjes integral, using
integration by parts and the smoothness of the Poisson kernel.
Recall that the Poisson kernel
$$
\Poi(z,w)= \frac{1}{2\pi}\Re \left(\frac{1+z\overline w}{1-
z\overline w} \right)= \frac{1}{2\pi}\Re\left( \frac{2}{1-
z\overline w} -1 \right) = 2\Re S_\UU(z,w)-\frac{1}{2\pi}
$$
has the kernel property
$$
\Poi(z,w)=\intc \Poi(z,e^{it})\Poi(e^{it},w) \,dt \,.
$$
(This follows  from the Poisson formula for harmonic
functions, see \cite{ahlfors}, Section 6.3). The white
noise $dB$ has the property that if $f_1,\ f_2$ are
smooth functions on an interval and $ f^\#_i = \int
f_i(t) \,d B(t) $ then $ \ev f_1^\#f_2^\# =  \int
f_1(t) f_2(t) \,dt. $ By this and the kernel property
we get $ \ev \Bigl(u(z)u(w) \Bigr) = \Poi(z,w). $
Therefore if $b$  is a standard real Gaussian
independent of $B(\cdot)$, then \bel{bterm}
 \tilde u(z)=\sqrt{\frac{\pi}{2}} u(z)+\frac{b}{2}
\ee
 has covariance structure
 $
\ev [\tilde u(z) \tilde u(w)] = \pi \Re\,S_\UU(z,w).
 $
Now if $\nu,\ \nu'$ are mean $0$ complex
Gaussians, then $\Re \ev \nu\overline{\nu'} = 2\ev (\Re \nu \Re
\nu')$; thus (\ref{covariance}) implies  that
 $\tilde u$ has the same distribution as $\Re f_\UU$.

\medskip
\noindent{\bf Remark.}\  Similarly, since $f_{\UU,2}$
is the derivative of $\sum_{m=1}^\infty a_m
z^m/{\sqrt{m}}$, the zero set $Z_{\UU,2}$ can be
interpreted as the set of saddle points of the random
harmonic function
 $$u(z)=\sum_{m=1}^\infty \Re(a_m z^m)/{\sqrt{m}}$$
in $\UU$. More generally, in any domain $D$, the zero
set $Z_{D,2}$ can be interpreted as the set of saddle
points of the Gaussian free field (with free boundary
conditions) restricted to harmonic functions.

\medskip

\noindent{\bf Joint moments of complex Gaussians.} We
will need the following known fact for the proofs of
Hammersley's formula and Theorem \ref{detformula}.

\begin{fact}\lb{permanent}
If $Z_1,\ldots,Z_n$ are mean 0 jointly complex
Gaussian random variables with covariance matrix
${\mathcal M}_{jk}=\ev Z_j \overline{Z_k}$, then
 $\ev \big(|Z_1\cdots Z_n|^2\big) = \perm(\mathcal M)$.
\end{fact}
\begin{proof}
We will check that in general for jointly complex
normal, mean 0 random variables $Z_j$, $W_j$ we have
$$
\ev[Z_1 \cdots Z_n \overline W_1 \cdots \overline W_n]
= \sum_{\sigma} \prod_{j=1}^k \ev Z_j \overline
W_{\sigma(j)} = \perm \left(\ev Z_j\overline
W_k\right)_{jk},
$$
where the sum is over all permutations $\sigma\in
S_n$. (See the book of \cite{simon} for a similar
statement in the real case.) Both sides are linear in
each $Z_j$ and $\overline W_j$, and we may assume that
the $Z_j$, $W_j$ are complex linear combinations of
some finite i.i.d.\ standard complex Gaussian sequence
$\{V_j\}$. The formula is proved by induction on the
total number of nonzero coefficients that appear in
the expression of the $Z_j$ and $W_j$ in terms of the
$V_j$. If the number of nonzero coefficients  is more
than one for one of $Z_j$ or $W_j$, then we may write
that variable as a sum and use induction and
linearity. If it is 1 or 0 for all $Z_j,\,W_j$, then
the formula is straightforward to verify; in fact,
using independence it suffices to check that $V=V_j$
has  $\ev V^n \overline V^m = n!\one_{\{m=n\}}$. For
$n\not=m$ this follows from the fact that $V$ has a
rotationally symmetric distribution. Otherwise,
$|V|^{2n}$ has the distribution of the $n$th power of
a rate 1 exponential random variable, so its
expectation equals $n!$.
\end{proof}

\section{A determinant formula in the i.i.d.\ case}

The goal of this section is to prove Theorem
\ref{detformula} and Corollary \ref{cdetformula}. The
proof relies on the  i.i.d.\ nature of the
coefficients of $f=f_{\UU}$, M\"obius invariance,
Hammersley's formula and Borchardt's identity
\re{borchardt}.

For $\beta \in \UU$ let
 \bel{moeb}
 T_\beta(z) = \frac{z -\beta }{1-
 \overline \beta z}
 \ee
 denote a M\"obius transformation fixing the unit
 disk, and define
$$
\tau_\beta(z)= \frac{(1 - |\beta|^2)^{1/2}}{1-\overline{\beta} z}
\,,
$$
so that $\tau_\beta^2(z)=T_\beta'(z)$.
\begin{remark}\lb{projections} Recall that for two jointly
 complex Gaussian random
vectors $X,Y$, the distribution of $Y$ given $X=0$ is
the same as the distribution of $Y$ with each entry
projected to the orthocomplement
(in $L^2$ of the underlying probability space)
of the subspace spanned by the components $X_i$ of $X$.
\end{remark}

\begin{proposition}\label{essential} Let
$f=f_\UU$ and $z_1,\ldots,z_n \in \UU$.
The distribution of the random function
 \bel{essence}
 T_{z_1}(z)\cdots T_{z_n}(z)f(z)
 \ee
is the same as the conditional distribution of $f(z)$ given
$f(z_1)=\ldots =f(z_n)=0$.
\end{proposition}

\begin{proof}
First consider $n=1$.
The assertion is clear for $z_1=0$;
here the i.i.d.\ property of the $\{a_k\}$ is crucial.
More generally, for $z_1=\beta$,  by
(\ref{gentran}) the random function $\wf=\tau_\beta \cdot (f\circ T_{\beta})$
has the same distribution as $f$. Since $T_\beta(\beta)=0$, from the formula
$$
\wf(z)= \tau_\beta(z)  \sum_{k=0}^\infty a_k
\bigl(T_\beta(z)\bigr)^k
$$
it is clear that the distribution of $T_\beta \cdot \wf$
is identical to the conditional distribution
of $\wf$  given  $\wf(\beta)=0$, whence the same must hold for $f$
in place of $\wf$.
The proposition for $n > 1$ follows by induction: to go from $n$ to $n+1$,
we must condition $(f \gi f(z_1)=\ldots=f(z_n)=0)$ on
$f(z_{n+1})=0$. By the assumed identity for $n$ points, this is
 equivalent to conditioning $(T_{z_1} \cdots T_{z_n} \cdot f)(z)$ on
$f(z_{n+1})=0$. By Remark \ref{projections},
conditioning is a linear operator that commutes with multiplication by
the deterministic functions $T_{z_i}$. Applying the
equality of distributions $(f(z) \gi f(z_{n+1})=0) \ed T_{z_{n+1}}(z) f(z)$
completes the proof.
\end{proof}

Fix $z_1,\ldots,z_n \in \UU$ and denote
\bel{defups}
\Upsilon(z)= \prod_{j=1}^n T_{z_j}(z) \,.
\ee
Since $T_{z_k}(z_k)=0$ and $T'_{z_k}(z_k)=1/(1-z_k\zb_k)$, we have
\bel{derups}
\Upsilon'(z_k)=T_{z_k}'(z_k) \cdot \prod_{ \, j: \, j\not=k} T_{z_j}(z_k)=
\prod_{j=1}^n \frac{1}{1-z_j\zb_k} \cdot \prod_{ \, j: \, j\not=k} (z_j-z_k)
\ee
for each $k \le n$.

\begin{corollary}\lb{dercor}
Let $f=f_\UU$ and $z_1,\ldots,z_n \in \UU$. The conditional joint
distribution of the random variables $\big(f'(z_k) \, : \,
k=1,\ldots, n \big) $
 given that $f(z_1)=\ldots =f(z_n)=0$,
is the same as the unconditional joint distribution of \
$\big(\Upsilon'(z_k)f(z_k)\, : \, k=1,\ldots, n \big) \,. $
\end{corollary}

\begin{proof}  The conditional distribution of $f$
given that $f(z_j)=0$ for $1 \le j \le n$, is the same as the
unconditional distribution of $\Upsilon \cdot f$.
Since $\Upsilon(z_k)=0$, the derivative of
$\Upsilon(z)f(z)$ at $z=z_k$ equals $\Upsilon'(z_k)f(z_k)$.
\end{proof}

Consider the $n \times n$ matrices $A$ and $M$,
 with entries
\begin{eqnarray*}
\label{jdensity} A_{jk}&=& \ev f(z_j)\overline{f(z_k)}  = (1-z_j\overline z_k)^{-1},\\
M_{jk}&=&(1-z_j\overline z_k)^{-2}.
\end{eqnarray*}
By the classical Cauchy determinant
formula, see \cite{muir} p.\ 311, we have
 $$
 \det(A)=\prod_{k,j}\,\frac{1}{1-z_j \zb_k} \;
 \prod_{k<j}\, (z_j-z_k)(\zb_j-\zb_k) \,.
 $$
Comparing this to \re{derups}, we see that
\bel{cauchv} \det(A)=\prod_{k=1}^n |\Upsilon'(z_k)|
\,. \ee We will need the classical identity of
\cite{borchardt1855} (see also \cite{minc}):
\bel{borchardt}
 \perm\left(\frac{1}{x_j+y_k}\right)_{j,k} \cdot
\det\left(\frac{1}{x_j+y_k}\right)_{j,k} =
 \det \left(\frac{1}{(x_j+y_k)^2}\right)_{j,k} \,.
 \ee
Setting $x_j=z_j^{-1}$ and $y_k=-\overline{z_k}$ and
dividing both sides by $\prod_j z_j^2$, gives that
\bel{borch2}
 \perm(A)\det(A)=\det(M).
\ee
We are finally ready to prove Theorem \ref{detformula}.
Corollary \ref{cdetformula} is a direct consequence of
the conformal invariance in Proposition \ref{ci} and
the way the \szego\ and Bergman kernels transform  under conformal
maps (see \re{szegoci}).

\begin{proofof}{Theorem \ref{detformula}}
Recall from \re{prehan} that \bel{suf1} p(z_1,\ldots,z_n)=
\frac{\ev\Bigl(|f'(z_1)\cdots f'(z_n)|^2 \, \Big| \,
f(z_1),\ldots, f(z_n)=0 \Bigr)}{\pi^n \det (A)} \,.
 \ee
By Corollary \ref{dercor}, the numerator of the right hand side of
\re{suf1} equals
 \bel{suf2}
  \ev\big(|f(z_1)\cdots
 f(z_n)|^2\big) \,\prod_{k}
|\Upsilon'(z_k)|^2 = \perm(A) \det(A)^2 \,,
 \ee
where the last equality uses the
Gaussian moment formula of Fact \ref{permanent}
and \re{cauchv}. Thus
$$
p(z_1,\ldots,z_n)=\pi^{-n} \perm(A)\det(A)=\pi^{-n} \det(M),
$$
by \re{borch2}.
\end{proofof}

\section{The number of zeros of $f_\UU$ in a disk} \lb{s.number}

In this section we prove Theorem \ref{distribution}
and Corollary \ref{hole}. In fact, the corollary only
uses part (i) of the theorem, so we delay the proof of
part (ii) to the end of the section.

\begin{lemma}\label{funmom}
Let $r_1\le\ldots \le r_m$, let $B_j=B_{r_j}(0)$, and
let $\wN_j=\#(Z_{\UU}\cap B_j)$. Then
 \bel{cumu}
 \ev \wN_1(\wN_2-1) \cdots (\wN_m-m+1) = \sum_{\sigma} \prod_{\nu\in \sigma}
(-1)^{|\nu|+1}\frac{r_\nu^2}{1-r_\nu^2},
 \ee
where the sum is over all permutations $\sigma$ of
$\{1,\ldots ,m\}$, the product is over all cycles
$\nu$ of the permutation, $|\nu|$ is the length of $\nu$,
and $r_\nu=\prod_{j\in \nu} r_{j}$.
\end{lemma}

\newcommand{\dzms}{\,dz_1\cdots dz_m}
\begin{proof}
Applying Proposition \ref{permformula} (ii) to the set
$B_1\times \cdots\times B_n$ we get
\begin{eqnarray*}
\ev \wN_1(\wN_2-1) \cdots (\wN_m-m+1) &=& \int_{B_1\times
\cdots \times B_m}
 p(z_1,\ldots,z_{m}) \dzms \\
 &=& \int_{B_1\times \cdots \times B_m}\det \big(K(z_i,z_j)\big)_{i,j}
 \dzms.
\end{eqnarray*}
Expanding the determinant and exchanging sums and
integrals we get a sum over all permutations of $m$
elements:
$$
\sum_\sigma \operatorname{sgn}(\sigma)\int_{B_1\times
\cdots \times B_m} K(z_1,z_{\sigma_1}) \cdots
K(z_m,z_{\sigma_m}) \dzms.
$$
For each permutation $\sigma$, the corresponding
integral is a product over cycles $\nu$ of $\sigma$ of
 \bel{kcycle}
I_\nu=\int K(z_1,z_2)K(z_2,z_3)  \cdots
K(z_{|\nu|},z_1) \,dz_1 \cdots dz_{|\nu|}.
 \ee
where $\nu$ is an ordered subset of $\{1,\ldots, m\}$
and each variable $z_i$ ranges over the disk
$B_{\nu_i}$ of radius $r_{\nu_i}$.  The formula for
the Bergman kernel gives
$$
K(z_1,z_2)=\frac{1}{\pi(1-z_1\overline
z_2)^{2}}=\pi^{-1}\sum_{n=0}^\infty (n+1)(z_1\overline
z_2)^n.
$$
Using this, we expand the product in \re{kcycle} into
a sum of monomials in the variables  $\{z_j\}_{j=1}^n$
and $\{\overline z_j\}_{j=1}^n$; then we
 integrate term by term. Each monomial  in
which the exponents of $z_j$ and $\overline z_j$ are
different for some $j$ integrates to 0. Thus in all
remaining terms the exponents of $z_j$ and $\overline
z_j$ are the same. Since $z_j$ always comes as part of
a product $z_j\overline z_{j+1}$, the exponents of
$z_j$ and $\overline z_{j+1}$ have to be the same as
well. This implies that in nonvanishing terms all
exponents agree, and we are left with
\begin{eqnarray*}
 I_\nu &=& \sum_{n=0}^\infty
 \pi^{-|\nu|}\int (n+1)^{|\nu|} \big(z_1\overline z_1 \cdots
 z_{|\nu|}\overline z_{|\nu|}\big)^n\,dz_1 \cdots dz_{|\nu|} \, .
\end{eqnarray*}
Since
$$ (n+1)\int_{B_r(0)}
|z|^{2n}\,dz = 2\pi(n+1)\int_0^r s^{2n+1}\,ds =\pi
r^{2n+2},
$$
setting $r_\nu=\prod_{j\in \nu} r_j$ we get
$$
I_\nu = \sum_{n=0}^\infty r_\nu^{2(n+1)}
=\frac{r_\nu^2}{1-r_\nu^2}.
$$
Since $\operatorname{sgn}(\sigma)=\prod_{\nu \in\sigma} (-1)^{|\nu|+1}$,
this completes the proof of the lemma.
\end{proof}

\begin{proofof}{Theorem \ref{distribution} (\rm i)}
Denote $\bb_k =\ev \binom{N_r}{k}$, and
 let $\mathcal P$ be chosen uniformly
at random from all permutations of $\{1,\ldots ,k\}$.
Let $q=r^2$. Then by \re{cumu}, we have
$$
\bb_k =\ev \prod_{y\in \mathcal P}
(-1)^{|y|+1}
\frac{q^{|y|}}{1-q^{|y|}} \,,
$$
where  the product is over  cycles $y$ of
 $\mathcal P$.  Since the cycle containing $1$ of
$\mathcal P$ has length that is uniform on $\{1,\ldots
,k\}$, and given that cycle, the other cycles form a
uniform permutation on the rest of $\{1,\ldots ,k\}$, we
get the recursion
\bel{bbform}
\bb_k = \frac{1}{k}\sum_{\ell=1}^k
(-1)^{\ell+1}\frac{q^\ell}{1-q^\ell}\, \bb_{k-\ell} \,,
\ee
with $\bb_0=1$.
Consider the generating function $\bb(s)=\sum_{k \ge
0}\, \bb_k s^k$.
Multiplying \re{bbform} by $ks^k$ and
summing over $k\ge 1$, we get
\bel{diffe}
s\bb'(s)=\bb(s)(s\psi(s)) \,,
\ee
where
$$
\psi(s)=\sum_{\ell=1}^\infty
(-s)^{\ell-1}\frac{q^{\ell}}{1-q^{\ell}}=
\sum_{k=1}^\infty \sum_{\ell=1}^\infty(-s)^{\ell-1}
q^{k\ell}.
$$
We write \re{diffe} in the form
$$
(\log \bb(s))'=\psi(s),
$$
which we integrate to get
$$
\log \bb(s)=
-\sum_{k=1}^\infty \sum_{\ell=1}^\infty \frac{(-q^ks)^\ell}{\ell}
=\sum_{k=1}^\infty \log(1+q^ks) \, ,
$$
where the constant term is zero as $\bb_0=\bb(0)=1$. Thus
\bel{prepeuler}
\bb(s)=\prod_{k=1}^\infty (1+q^ks) \,.
\ee

Taking expectations of the identity
$$
(s+1)^{N_r} = \sum_{k=0}^\infty \binom{N_r}{k} s^k
$$
gives
$$
\ev(s+1)^{N_r} = \sum_{k=0}^\infty \bb_k s^k
=\bb(s) \, ,
$$
and this concludes the proof of (i).
\end{proofof}

%

\begin{proofof} {Corollary \ref{hole}}
\noindent{\bf (i)}
Theorem \ref{distribution} implies that
$\pr(N_r=0)=\prod_{k=1}^\infty (1-r^{2k})$
and the asymptotics for the right hand side are classical,
see \cite{newman}, p.\ 19.
For the reader's  convenience we indicate the argument.
Let $
L=\log \pr(N_r=0)=\sum_{k=1}^\infty \log(1-r^{2k})$
which we compare to the integral
\bel{logexp}
I=\int_1^\infty \log(1-r^{2k})\,dk
=\frac{1}{-2\log r}\int_{-2\log r}^\infty
\log(1-e^{-x})\,dx.
\ee
We have $I+\log(1-r^2)< L < I$, so $L=I+o(h)$. Since
$
-\log(1-e^{-x})=\sum_{n=1}^\infty \frac{e^{-nx}}{n},
$
the integral in \re{logexp} converges to $-\pi^2/6$.
But
$\frac{-1}{2\log r}=\frac{1/2+o(1)}{1-r}=\frac{h}{4\pi}+o(h),$
and we get
$
L=-\frac{\pi^2/12+o(1)}{1-r}=-\frac{\pi h}{24} +o(h) \,,$ as claimed.

\smallskip

\noindent{\bf (ii)} One of Euler's partition
identities (see \cite{paksurvey}, section 2.3.4) gives
\bel{euler} \prod_{k=1}^\infty (1+q^ks) =
\sum_{k=0}^\infty
\frac{q^{\binom{k+1}{2}}s^k}{(1-q)\cdots (1-q^k)}. \ee
so the claim follows from \re{prepeuler}.
\smallskip

\noindent{\bf (iii)} The formulas for the mean and variance of
$N_r$ follow from the binomial moment formula \re{cumu}. Using the
general central limit theorem due to \cite{cl95} and
\cite{soshnikov00b} (p.\ 497) for determinantal processes, we get
that as $r\to 1$, the normalized distribution of $N_r$ converges
to standard normal, as required. Alternatively, the last claim can
be easily verified by computing the asymptotics of the moment
generating function directly. Yet another way is to apply
Lindeberg's triangular array central limit theorem to the
representation  of $N_r$ as the sum of independent random
variables, as given in Theorem \ref{distribution}(i).
\end{proofof}

\subsection*{The joint distribution of the moduli of
$Z_\UU$}

{\noindent{\bf Proof of part (ii) of Theorem \ref{distribution}:}

 The zero set of $f_\UU$ is determinantal with the Bergman kernel
$K(z,w)$.
Let $$
K_n(z,w)=\frac{1}{\pi}\sum_{j=0}^{n-1} (j+1)(z\overline w)^j
\,.  $$
Since $K_n(z,w) \to K(z,w)$ as $n \to \infty$ uniformly on
compact sets of $\UU^2$, Proposition 3.10 of \cite{ST} yields that the
determinantal point processes with kernels $K_n$ converge weakly, as
$n \to \infty$, to $Z_\UU$.  Thus it suffices to prove that the set of
absolute values $\{|\zeta_j|\}_{j=1}^n$ of the $n$ random points of
the determinantal process with kernel $K_n$, has the same law as
$\{U_j^{\frac{1}{2j}}\}_{j=1}^n$ where $U_j$ are i.i.d.\ uniform on
[0,1].

For any $z_1,\ldots ,z_n$,

\[ \l( \begin{array}{cccc}
      K_n(z_1,z_1) & \ldots &\ldots & K_n(z_1,z_n) \\
      \ldots       & \ldots &\ldots & \ldots       \\
      \ldots       & \ldots &\ldots &\ldots       \\
      K_n(z_n,z_1) & \ldots &\ldots &K_n(z_n,z_n) \end{array} \r)= \]
\[
 \frac{1}{\pi^n} \l( \begin{array}{cccc}
          1  & z_1 & \ldots & z_1^{n-1} \\
          1  & z_2 & \ldots & z_2^{n-1} \\
      \cdots & \cdots  & \cdots & \cdots \\
          1  & z_n & \ldots & z_n^{n-1}  \end{array} \r)
 \l( \begin{array}{cccc}
          1  & 0 & \ldots & 0 \\
          0 & 2 & \ldots & 0 \\
      \cdots & \cdots  & \ddots & \cdots \\
          0 & 0 & \ldots & n  \end{array} \r)
 \l( \begin{array}{ccc}
         1          & \ldots & 1 \\
       \overline z_1  & \ldots & \overline z_n \\
      \cdots        & \cdots & \cdots  \\
     \overline z_1^{n-1}&\ldots &\overline z_n^{n-1} \end{array} \r).
\]
Setting $z_j=r_je^{i\theta_j}$ we find that the joint intensity
of $\{|\zeta_j|\}_{j=1}^n$, evaluated at $\{r_j\}_{j=1}^n$, equals
\begin{eqnarray} \lb{newcorr}
& & \intt_{[0,2\pi]^n} \det\Bigl(K_n(z_j,z_k) \Bigr)_{j,k=1}^n
 \, r_1\,d\theta_1\cdots r_n\,d\theta_n
 \\[1.5ex] \nonumber
& = & \frac{n!}{\pi^n} \intt \Bigl(\sum_{\sig} \mb{sgn}(\sig) \prodd_{j=1}^n z_j^{\sig_j-1}
\Bigr)
  \Bigl(\sum_{\tau} \mb{sgn}(\tau) \prodd_{j=1}^n \overline
  z_j^{\tau_j-1} \Bigr) \, r_1\,d\theta_1\cdots r_n\,d\theta_n\,.
\end{eqnarray}
When we expand the sums, for $\sig\not=\tau$ the
integrand contains a factor of the form $z_j^p\overline z_j^q$ with
$p\not=q$, and therefore the integral vanishes. When $\sig=\tau$, we get
$(2\pi)^{n}\prodd_j r_j^{2\sig_j-1}$. Thus \eqref{newcorr} equals
\begin{equation}\label{eq:randomizedpdf} 2^n
  n!\summ_{\sig}\prodd_{j=1}^n r_j^{2\sig_j-1}. \end{equation}
Now $U_j^{\frac{1}{2j}}$ has density $2jx^{2j-1}$ in $[0,1]$. Hence, the
joint intensity
of $\{U_1^{\frac{1}{2}},\ldots , U_n^{\frac{1}{2n}}\}$
is precisely  (\ref{eq:randomizedpdf}).   This proves the theorem.

\noindent{\bf Remark.} The proof above is modeled after an argument of
\cite{kostlan}
for the distribution of the eigenvalues of a random complex Gaussian matrix.
It is simpler than our original proof that relied on random permutations.

\section{Law of large numbers} \lb{s.lln}

The goal of this section is to prove Proposition \ref{lln},
the law of
large numbers for the number of zeros of $f_{\UU,\rho}$.
We will use the following lemma in the proof.

\begin{lemma}\lb{ballbest}
Let $\mu$ be a Borel measure on
a metric space $S$, and
assume that all balls of the same radius have the same measure.
Let $\psi:[0,\infty) \to [0,\infty)$ be a non-increasing function.
 Let $A \subset S$ be a Borel set, and let $B=B_R(x)$ be a ball
centered at $x\in S$ with $\mu(A)=\mu(B_R(x))$. Then for all $y \in S$
$$
\int_A \psi(\dist(y,z))\,d\mu(z) \le \int_B \psi(\dist(x,z))
\,d\mu(z).
$$
\end{lemma}
\begin{proof}
It suffices to check this claim for indicator functions
$\psi(x)=\one_{\{x\le r\}}$. In this case, the inequality reduces to
$$
\mu(A\cap B_r(y)) \le \mu(B_R(x)\cap B_r(x)),
$$
which is clearly true  both for $r\le R$ and for
$r>R$.
\end{proof}

\begin{proofof}{Proposition \ref{lln}}
We have
$$
\ev N(\h) = \int_\Lambda p(z)\,dz = \frac{4\rho}{\pi} h \, .
$$
Let $Q(z,w)=p(z,w)/(p(z)p(w))$.
Then by formula \re{2point} we have
$$
Q(0,w) -1 \le C(1-|w|^2)^\rho \,.
$$
we denote the right hand side by $\psi(0,w)$ and extend $\psi$ to
$\UU^2$ so that it only depends on hyperbolic distance.
\begin{eqnarray*}
 \ev(N(\h)(N(\h)-1))-(\ev N(\h))^2 &=&
 \int_\Lambda \int_\Lambda
   \big(p(z,w)-p(z)p(w)\big)\,dw\,dz \\ &=&
   \int_\Lambda \int_\Lambda (Q(z,w)-1) \,p(w)\,dw\, p(z)\,dz  \\ &\le&
   \int_\Lambda \int_\Lambda \psi(z,w)\,p(w)\,dw\, p(z)\,dz
\end{eqnarray*}
Let $B_R(0)$ be a ball with hyperbolic area $h=4\pi R^2/(1-R^2)$. Note
that $p(w)dw$ is constant times the hyperbolic area element, so we may
use Lemma \ref{ballbest} to bound the inner integral by
\begin{eqnarray*}
\int_{B_R(0)} \psi(0,w)\,p(w)\,dw &=& c \int_0^R
(1-r^2)^\rho (1-r^2)^{-2} r \,dr\\
&=& \frac{c}{2}\int_S^1 s^{\rho-2} \,ds
\end{eqnarray*}
with $S=1-R^2$. Thus we get
 \bel{42}
 \Var\ N(\h) = \ev N(\h)+\ev(N(\h)(N(\h)-1))-(\ev N(\h))^2
 \le \frac{h\rho}{4\pi}+
 \frac{ch\rho}{8\pi}\int_S^1 s^{\rho-2}\,
 ds.
 \ee
For $\rho>1$ this is integrable, so $\Var\ N(\h) \le O(h)$. For
$\rho<1$ we can bound the right hand side of \re{42} by
$O(hS^{\rho-1})=O(h^{2-\rho})$. Thus in both cases, as well as
when $\rho=1$ (see Corollary \ref{hole}(iii)), we have
$$ \Var\ N(\h) \le c (\ev N(\h))^{2-\beta}
$$
with $\beta=\rho \wedge 1>0$. For $\eta>1/\beta$, we find that
$$
Y_k=\frac{N(k^\eta)-\ev N(k^\eta)}{\ev N(k^\eta)}
$$
satisfies  $\ev Y_k^2=O(k^{-\eta\beta})$, whence $ \ev \sum_k
Y_k^2<\infty, $ so $Y_k\to 0$ a.s. Monotonicity and interpolation
now give the desired result.
\end{proofof}

\section{Reconstructing the function from its zeros}
\lb{s.rec}
The goal of this section is to prove Theorem \ref{recprod}.
The main step in the proof is the following.
\begin{proposition}\lb{jensen} Let
$c'_\rho =e^{\rho/2-\gamma/2}$. We have
$$
|f_{\UU,\rho}(0)|=c'_\rho \lim_{r\to 1}
\;(1-r^2)^{-\rho/2}\hspace{-2pt}\prod_{\substack{z\in Z_{\UU,\rho}
\\ |z|<r}} \hspace{-2pt}|z| \as
$$
\end{proposition}

We first need a simple lemma.
\begin{lemma}\lb{covariancelemma}
If $X$, $Y$ are jointly complex Gaussian with variance
1, then for some absolute constant $c$ we have
 \bel{coveq}
 \Bigl|\Cov\Bigl(\hspace{-2pt}\log|X| \, ,\, \log|Y|\Bigr)\Bigr|\le c
\Bigl|\ev(X\overline Y)\Bigr|.
 \ee
\end{lemma}

\begin{proof}
Write $Y=\alpha X + \beta Z$, where $X,\,Z$ are
i.i.d.\ standard complex Gaussian variables,
$\alpha=E(X\overline Y)$ and $|\alpha|^2+|\beta|^2=1$.
It clearly suffices to consider $|\alpha|<1/2$. Since
$$
\log|Y|=\log|\beta Z|+\log |1+\alpha X/(\beta Z)|,
$$
 the inequality \re{coveq} reduces to
 \bel{coveq2}
 \Bigl|\Cov \Bigl(\hspace{-2pt}\log|X| \, , \,
 \log\Bigl|1+\frac{\alpha X}{\beta Z}\Bigr|\Bigr)\Bigr|\le c|\alpha|.
 \ee
We will use the estimate
 \bel{coveq3}
 \ev\big|\hspace{-2pt}\log|1+\lambda/Z|\big|\le
 c_1|\lambda| \mbox{\ \ \ for } \lambda \in
 \mathbb{C},
 \ee
which can be verified by considering  the positive and
negative parts of $ \log|1+\lambda/Z|$ as follows. The
positive part is handled using the numerical
inequality $\log|1+w| \le |w|$ and the integrability
of $|Z|^{-1}$. For the negative part, when
$|\lambda|\ge 1$, the density of $|1+\lambda/Z|$ is
uniformly bounded in the disk of radius $1/2$, so it
remains to consider the case $|\lambda|< 1$. Then $
\ev \log_{-}|1+\lambda/Z|$
 can be controlled by partitioning into the events
 $$
 G_k=\{e^{-k}< |1+\lambda/Z| \le e^{1-k} \} \,.
 $$
 Since $\pr(G_k)=O(|\lambda^2| e^{-2k})$, we get
 $$
  \ev \Bigl(\one_{G_k} \log_{-}|1+\lambda/Z|  \Bigr)
 = O(k |\lambda^2| e^{-2k}) \,.
 $$
 Summing over $k$ establishes \re{coveq3}.

By conditioning on $X$, \re{coveq3} yields
$$
\ev\Bigl(\log |X| \cdot \log \Bigl|1+\frac{\alpha X}{\beta Z}\Bigr| \Bigr)
 \le c_1
\Bigl|\frac{\alpha}{\beta}\Bigr| \cdot \ev \big|X\log|X|\big|=c_2
\Bigl|\frac{\alpha}{\beta}\Bigr|.
 $$
This bounds the first term (expectation of the product) in the covariance
on the left hand side of
\re{coveq2}. The second term (product of expectations)
can be bounded by the same argument.
\end{proof}

\begin{proofof}{Proposition \ref{jensen}}
 Assume that $f=f_{\UU,\rho}$ has no zeros at $0$ or on the
circle of radius $r$. Then Jensen's formula
(\cite{ahlfors}, Section 5.3.1) gives
$$
\log |f(0)|=\frac{1}{2\pi} \int_0^{2\pi} \log
|f(re^{i\alpha})|\;d\alpha +\sum_{z\in Z,\, |z|<r}\log
\frac{|z|}{r} \,,
$$
where $Z=Z_{\UU,\rho}$. Let
$|f(re^{i\alpha})|^2=\sigma_r^2 Y$, where
$$\sigma_r^2=\Var
f(re^{i\alpha})=[2\pi
S_\UU(r,r)]^\rho=(1-r^2)^{-\rho}$$ and $Y$ is an
exponential random variable with mean 1. We have
$$\ev \log |f(re^{i\alpha})|=\frac{\log \sigma_r^2+ \ev \log Y}{2}=
\frac{-\rho\log(1-r^2)-\gamma}{2},
$$
where the second equality follows from the integral
formula for Euler's constant
$$
\gamma=-\int_0^\infty e^{-x}\log x \,dx.
$$
Introduce the notation
$$
g_r(\alpha)=\log |f(e^{i\alpha}r)| + \frac{\rho\log
(1-r^2)+\gamma}{2}
$$
so that the distribution of $g_r(\alpha)$ does not
depend on $r$ and $\alpha$, and $\ev g_r(\alpha)=0$.
Let
$$
L_r= \frac{1}{2\pi} \int_0^{2\pi} g_r(\alpha) d\alpha.
$$
We first prove that $L_r\to 0$ a.s. over a suitable
deterministic sequence $r_n \uparrow 1$. We compute:
$$
\Var\ L_r=\ev \left( \frac{1}{(2\pi)^2}
\int_0^{2\pi}\int_0^{2\pi}
g_r(\alpha)g_r(\beta)\,d\beta\,d\alpha\right).
$$
Since the above is absolutely integrable, we can
exchange integral and expected value to get
$$
\Var\ L_r= \frac{1}{(2\pi)^2}
\int_0^{2\pi}\int_0^{2\pi}
\ev(g_r(\alpha)g_r(\beta))\,d\beta\,d\alpha =
\frac{1}{2\pi} \int_0^{2\pi}
\ev(g_r(\alpha)g_r(0))d\alpha.
$$
where the second equality follows from rotational
invariance. By Lemma \ref{covariancelemma}, we have
$$
\ev \big( g_r(\alpha)g_r(0) \big)\le c
 \frac{\bigl|\ev \bigl(f(re^{i\alpha})
  \overline{f(r)}\bigr)\Bigr|}{\Var\,f(r)}=
c\Bigl|\frac{1-r^2}{1-r^2e^{i\alpha}}\Bigr|^\rho \,.
$$
Let $\eps=1-r^2 <1/2$. Then for $\alpha\in[0,\pi]$ we can
bound
$$
{|1-r^2e^{i\alpha}|} \ge
 \left\{\begin{array}{ll}
 \eps &
 |\alpha| \le \eps \\
 2r^2 \sin \frac{\alpha}{2} \ge \frac{\alpha}{2} &
 \eps<\alpha<\pi/2 \\
 1&
 \pi/2 \le\alpha \le \pi,
\end{array}\right.
$$
which gives
$$
\frac{1}{c2\eps^{\rho}}\Var\ L_r \le \int_0^\pi
\frac{d\alpha}{|1-r^2e^{i\alpha}|^\rho}
 \le \eps^{1-\rho} +
 \frac{1}{2}\int_\eps^{\pi/2}\alpha\,d\alpha +
 \pi/2 \le
 \left\{\begin{array}{ll}
 c'&
 \rho<1 \\
 c'|\log \eps|
 &\rho=1 \\
 c'\eps^{1-\rho} &
 \rho>1.
\end{array}\right.
$$
By Chebyshev's inequality and the Borel-Cantelli
lemma,  this shows that, as $r \to 1$ over the
sequence $r_n=1-n^{-(1\vee( 1/\rho)+\delta)}$, we have
a.s.\ $L_{r_n} \to 0$ and
$$
\sum_{z\in Z,|z|<r}\log \frac{|z|}{r} -\frac{\rho\log
(1-r^2)+\gamma}{2} \to \log|f(0)|,
$$
or, exponentiating:
\bel{expon}
e^{-\gamma/2}(1-r^2)^{-\rho/2}
\hspace{-3pt}\prod_{\substack{z\in Z_{\UU,\rho} \\
|z|<r}} \hspace{-5pt}\frac{|z|}{r} \ \longrightarrow\
|f(0)|. \ee
 Since the product is monotone decreasing
and the ratio $(1-r_n^2)/(1-r_{n+1}^2)$ converges to
$1$, it follows that the limit is the same over every
sequence $r_n\to 1$ a.s.

Finally, by the law of large numbers (Proposition
\ref{lln}), the number of zeros $N_r$ in the ball of
Euclidean radius $r$ satisfies
 \bel{llnuse} N_r = \frac{r^2\rho}{1-r^2}\,(1+o(1))=
 \frac{\rho+o(1)}{1-r^2} \, \as,
 \ee
 whence
$$
r^{N_r} =\exp(N_r \log r)=e^{-\rho/2+o(1)} \, \as
$$
Multiplying this with \re{expon} yields the claim.
\end{proofof}

\begin{proofof}{Theorem \ref{recprod}} {\bf(i)}
By the law of large numbers for $N_r$ (see also \re{llnuse}),
 \bel{collect}
 \sum_{|z_k| \le r} \frac{1}{k}
 = \gamma + \log N_r +o(1)
 = \gamma+\log \rho -\log(1-r^2) +o(1) \,.
 \ee
 Multiplying by $\rho/2$ and exponentiating, we get that
\bel{collect2} \prod_{|z_k| \le r} e^{\rho/(2k)} =
e^{\gamma\rho/2}\rho^{\rho/2}(1-r^2)^{-\rho/2}(1+o(1))
\,. \ee In conjunction with Proposition \ref{jensen},
this yields \re{e.recp2}.

\medskip

\noindent{\bf(ii)} Let $f=f_{\UU,\rho}$ and
$$
T(z)=\frac{z-\zeta}{1-\overline{\zeta} z} \,.
$$
By \re{gentran}, $f$ has the same law as
\bel{samel}
\wf=(T')^{\rho/2} \cdot (f\circ T) \,.
\ee
Now $T'(\zeta)=(1-|\zeta|^2)^{-1}$. Therefore
$$
\wf(\zeta)=(1-|\zeta|^2)^{-\rho/2}f(0)= c_\rho
\prod_{k=1}^\infty e^{\rho/(2k)} |z_k| \quad \as,
$$
where $\{z_k\}$ are the zeros of $f$ in increasing modulus.
If $T(\zeta_k)=z_k$ then $\{\zeta_k\}$ are the zeros of $\wf$ in increasing
hyperbolic distance from $\zeta$. We conclude that
$$
\wf(\zeta)=c_\rho (1-|\zeta|^2)^{-\rho/2}
  \prod_{k=1}^\infty e^{\rho/(2k)} |T(\zeta_k)| \as
$$
\end{proofof}

For our study of the dynamics of zeros in the next
section, we will need a reconstruction formula for
$|f_{\UU,\rho}'(0)|$ when we condition on that $0 \in
Z_{\UU,\rho}$.

\begin{lemma} \lb{condz}
Denote by $\Omega_\eps$ the event that the power series
 $f_{\UU,\rho}$ defined in (\ref{genrho}) has a zero in $B_{\eps}(0)$.
As $\eps \to 0$, the conditional distribution of the
coefficients $a_1,a_2,a_3,\ldots$ given
$\Omega_\eps$, converges to a product law where $a_1$
is rotationally symmetric, $|a_1|$ has density
 $2r^3e^{-r^2}$, and $a_2,a_3,\ldots$
are standard complex Gaussian.
\end{lemma}

\begin{proof}

Let $a_0,\,a_1$ be i.i.d.\ standard complex normal
random variables, and $\rho>0$. Consider  first the limiting
distribution, as $\eps\to 0$, of $a_1$ given that the
equation $a_0+a_1 \sqrt{\rho}z=0$ has a root $Z$ in
 $B_\eps(0)$. The limiting
distribution must be rotationally symmetric, so it
suffices to compute its radial part. If $S=|a_0|^2$
and $T=|a_1|^2$, set $U=\rho|Z|^2= S/T$. The joint
density of $(S,T)$ is $e^{-s-t}$, so the joint density
of $(U,T)$ is $ e^{-ut-t} t$. Thus as $\eps\to 0$, the
conditional density of $T$ given $U<\rho\eps^2$ converges
to the conditional density given $U=0$, that is
$te^{-t}$. This means that the conditional
distribution of $a_1$ is not normal, rather, its radial
part has density $2r^3e^{-r^2}$.

\medskip
We can now prove the lemma.
The conditional density of the coefficients
$a_1,a_2,\ldots$ given $\Omega_\eps$, with respect to their original
product law, is given by the ratio
$\pr(\Omega_\eps \mid
a_1,a_2,\ldots)/\pr(\Omega_\eps)$. By Lemma
\ref{diffzero}, the limit of this ratio is not
affected if we replace $f_{\UU,\rho}$ by its linearization
$a_0+a_1 \sqrt{\rho}z$. This yields the statement of the lemma.
\end{proof}

Kakutani's absolute continuity criterion (see
\cite{williams}, Theorem 14.17) applied to the
coefficients gives the following
\begin{lemma}\lb{kakutani}
The distributions of the random functions $f_{\UU,\rho}(z)$ and
$(f_{\UU,\rho}(z)-a_0)/z$ are mutually absolutely continuous.
\end{lemma}

\begin{remark}\lb{prime}
By Lemma \ref{condz}, conditioning on $0\in
Z_{\UU,\rho}$ amounts to setting $a_0=0$ and changing
the distribution of $a_1$ in an absolutely continuous
manner. Thus, by Lemma \ref{kakutani}, given $0\in
Z_{\UU,\rho}$ the distribution of the random function
$g(z)=f_{\UU,\rho}(z)/z$ is absolutely continuous with
respect to the distribution of the unconditioned
$f_{\UU,\rho}(z)$. Hence we may apply Theorem
\ref{recprod} to $g(z)$ and get that given $0\in
Z_{\UU,\rho}$, if we order the {\em other\/} zeros of
$f_{\UU,\rho}$ in  increasing absolute value
  as $\{z_k\}_{k=1}^\infty$, then
\bel{jensenprime} |f'_{\UU,\rho}(0)|=|g(0)|=c_\rho
\prod_{k=1}^\infty e^{\rho/(2k)} |z_k| \quad \as \ee
\end{remark}

\section{Dynamics of zeros}\lb{s.dyn}

In order to understand the point process of zeros of
$f_\UU$ it is useful to think of it as a stationary
distribution of a time-homogeneous Markov process.

Define the complex Ornstein-Uhlenbeck process
$$
 a(t):=e^{-t/2}W(e^t),\ \ \ \
 W(t):=\frac{B_1(t)+iB_2(t)}{ \sqrt{2}},
 $$
where $B_1,\ B_2$ are independent standard Brownian
motions, and $W(t)$ is complex Brownian motion scaled
so that $\ev W(1)\overline{ W(1)}=1$. The process
$\{a(t)\}$ is then stationary Markov with the standard
complex Gaussian as its stationary distribution. First
we consider the process
$$
\varphi_t(z) \varphi_t(z; D)= \sum_{n=0}^\infty W_n(t) \psi_n(z), \ \
\ t>0
$$
where $W_n$ are independent complex Brownian motions and
$\{\psi_n(z)\}_{n\ge 0}$ is an orthonormal basis for
$H^2(D)$.  With $t= e^\tau$ we get the
time-homogeneous process,
$$
f_\tau(z) = e^{-\tau/2}\varphi_{e^\tau}(z)=
\sum_{n=0}^\infty a_n(\tau) \psi_n(z).
$$
Then the entire process $\varphi_t(z)$ (and so
$f_\tau(z)$) is conformally covariant in the sense that
if $T :\Lambda \to D$ is a conformal homeomorphism, then
the process
$$
\big\{\,\sqrt{T'(z)} \varphi_t(T(z))\, \big\}_{t>0}
$$
has the same distribution as $\varphi_t(z;\Lambda),\
t>0$. For this, by continuity, it suffices to check
that the covariances agree. Indeed, for $s \le t$,
$$
\ev \varphi_s(z)\overline{\varphi_t(w)}= \ev
\varphi_s(z)\overline{\varphi_s(w)}
$$
so the problem is reduced to checking the equality of covariances for a fixed
time, which has already been done in Proposition
\ref{ci}.

It follows automatically that the process $\{Z_D(t)\}$
of zeros of $\varphi_t$ is conformally invariant. To
check that it is a Markov process, recall from Section
\ref{s.ci} that $Z_D(t)$ determines $\varphi_t$ up to
a multiplicative constant of modulus 1. It is easy to
check that $\varphi_t$ modulo such a constant is a
Markov process; it follows that $Z_D(t)$ is a Markov
process as well.

\begin{remark}
This argument works in the case $\rho=1$. By replacing
the i.i.d.\ coefficients $a_n$ in \re{genrho} with OU
processes, it is possible to define a dynamic version
of the $\rho\not=1$ case in the unit disk. The same
argument as above shows that these are Markov
processes with distribution invariant under M\"obius
transformations of $\UU$.
\end{remark}

Finally, we give an SDE description of the motion of
zeros. Condition on starting at time $1$ with a zero
at the origin. This implies that $W_0(1)=0$, and by
the Markov property all the $W_i$ are complex Brownian
motions started from some initial distribution at time
1. For $t$ in a small time interval $(1,1+\eps)$ and
for $z$ in the neighborhood of $0$, we have
$$
\varphi_{t}(z)=W_0(t)+W_1(t) z+W_2(t) z^2 + O(z^3).
$$
If $W_1(1)W_2(1)\not=0$, then the movement the root
$z_t$ of $\varphi_t$ where $z_1=0$ is described by the
movement of the solution of the equation
$W_0(t)+W_1(t) z_t+W_2(t) z_t^2=O(z^3_t)$. Solving the
quadratic gives $$
 z_t=
 \frac{-W_1}{2W_2}
 \left(1-\sqrt{1-\frac{4W_0W_2}{W_1^2}}\right)+O(W_0^3).
$$
Expanding the square root we get
$$
z_t=-\frac{W_0}{ W_1} + \frac{W_0^2 W_2 }{W_1^3} +
O(W_0^3).
$$
Since $W_0(t)$ is complex, $W_0^2(t)$ is a martingale,
so there is no drift term. The noise term then has
coefficient $-1/W_1$, so the movement of the zero at
$0$ is described by the SDE (at t=1) $ dz_t=
-W_1(t)^{-1}dW_0(t) $ or, rescaling time for the
time-homogeneous version, for any $\tau$ with
$a_0(\tau)=0$ we get
 \bel{doesnot} dz_\tau = -\frac{1}{
a_1(\tau)} \,da_0(\tau).
 \ee
The absence of drift in \re{doesnot} can be understood
as follows: in the neighborhood we are interested in,
this solution $z_t$ will be an analytic function of
the $\{W_n\}$, and therefore has no drift.

For other values of $\rho$ the same argument gives
$$
dz_\tau = -\frac{1}{\sqrt{\rho} \,a_1(\tau)}\,
da_0(\tau).
$$

Of course, it is more informative to describe this
movement in terms of the relationship to other zeros,
as opposed to the coefficient $a_1$. For this, we
consider the reconstruction formula in Remark
\ref{prime}, which gives
$$|a_1|=|f'_{\UU,\rho}(0)|=c_\rho \prod_{k=1}^\infty
e^{\rho/(2k)} |z_k| \quad \as
$$
This means that when there are many other zeros close to a zero,
the noise term in its movement grows and it oscillates wildly.
This produces a repulsion effect for zeros that we have already
observed in the point process description. The equation
\re{doesnot} does not give a full description of the process as
the noise terms for different zeros are correlated.

\section{Hammersley's formula for Gaussian analytic functions}\lb{general}
%
A version of the following theorem was proved by
\cite{hammersley56}. The present version is from
\cite{friedman90}, Appendix B. We say that a point
process has {\bf integral joint intensity} $p$ if formula
\re{soshcor1} holds for its counting function $N$.

\begin{theorem} \lb{hammersley} Let
$f_n=a_nz^n+\ldots +a_0$ be a random polynomial, so
that $(a_0,\ldots, a_n)$ has an absolutely continuous
distribution with respect to Lebesgue measure on
$\CC^{n+1}$. Then the integral joint intensity of
zeros exists and  equals
\begin{equation}\lb{hamex}
p(z_1,\ldots,z_k)=\lim_{\eps\to 0} (\pi\eps^2)^{-k}
\int_{f(z_i)\in B_\eps(0),i=1,\ldots k}|f'(z_1)\cdots
f'(z_k)|^2\,da_0\cdots da_n.
\end{equation}
\end{theorem}

We also need the following consequence of Cauchy's
integral formula.

\begin{fact}\lb{cauchyfact}
Let $D$ be a bounded domain, and let $B\subset D$ be a
closed disk. Then for every $m\ge 0$ there exist
constants $c_m$ so that for every $f$ analytic on $D$
and every $z\in B$ the $m^{th}$ derivative satisfies
$|f^{(m)}(z)|\le c_m (\int_D |f(w)|^2 \,dw)^{1/2}$.
\end{fact}

\begin{proof} Cauchy's integral formula gives a uniform bound
on $f^{(m)}(z)$ for $z\in B$ in terms of the
$L^1$-norm of the function on any circle in $D$ about
$B$. Integration yields a bound in terms of the $L^1$
norm on an annulus, which is bounded above by the
$L^2$ norm on $D$.
\end{proof}

Next, we note some consequences of
the Taylor expansion for Gaussian analytic functions.

\begin{lemma}\lb{approxl}
Let $f$ be a Gaussian analytic function defined on a
domain $D$, and let $B\subset D$ be a closed disk
about $z_0$. Consider the partial sums of the Taylor
series expansion about $z_0$:
$$
f_n(z)=\sum_{k=0}^n a_k (z-z_0)^k.
$$
Then for all $m\ge 0$ the $m^{th}$ derivative satisfies
\begin{equation} \lb{sense}
\sup_B \ev |f^{(m)}_n-f^{(m)}|^2 \to 0 \; \mbox{ \rm as } n \to \infty \,.
\end{equation}
Consequently, for all $m_1,m_2\ge 0 $ the covariance function of
the derivatives of orders $m_1$ and $m_2$ of $f_n$ converges
uniformly on $B^2$ to the covariance function 
of the corresponding derivatives of $f$.
\end{lemma}

\begin{proof}
Note that finite a.s.\  limits of jointly Gaussian
random variables are are jointly Gaussian with finite
variance. This implies that the derivative of a
Gaussian analytic function $f$ is a Gaussian analytic
function. Moreover, the Taylor series of $f$ has
jointly Gaussian coefficients.
Consider the $L^2$ space of functions on the set
$X=\Omega \times D$ with the product of $\pr$ and
Lebesgue area measure.
Assume without loss of generality that $B$ is centered
at 0, and  let $f_n(z)=\sum_{j=0}^n a_jz^j$. Since
$f_n$ is a projection of $f$ in the space $X$ to the
subspace spanned by $f_0,\ldots, f_n$, it follows that
$f_n\to f$ in $L^2(X)$. By Fact \ref{cauchyfact}, we
have
$$
\ev \sup_B |f_n^{(m)}-f^{(m)}|^2 \to 0
$$
and therefore
$$
\sup_B \ev |f_n^{(m)}-f^{(m)}|^2 \to 0.
$$
which implies the claim for the covariance functions.
\end{proof}

\begin{corollary}\lb{approxcor}
Let  $g_n$ be polynomial of degree $n$ with i.i.d.\ standard
Gaussian coefficients independent of  $f_n$. Then $f_n+g_n/n!$
approximates $f$ in the sense of \eqref{sense} and for each
$n$ has coefficients with continuous joint density.
\end{corollary}

We first show a preliminary version of Proposition
\ref{permformula}; (i) will then follow from the
integral formula and a general lemma about point
processes.

\begin{proposition}\lb{prepermformula}
Using the notation of Proposition \ref{permformula},
assume that $A=A(z_1,\ldots, z_k)$ is nonsingular when
$z_1,\ldots,z_n\in D$ are distinct. Denote by
$N(\Lambda)$ the number of zeros of
 $f$ in $\Lambda$. 
 Then for any $n$ disjoint bounded Borel subsets
$\Lambda_1,\dots ,\Lambda_n$ of
 $D$, we have
\bel{soshcorpre} \ev \prod^n_{i=1} N(\Lambda_i)=
\int_{\Lambda_1\times\dots\times \Lambda_n}
 p(z_1,\dots ,z_n)\,dz_1\dots dz_n \,,
\ee where the integrations are with respect to
Lebesgue area measure.
\end{proposition}

\begin{proof}

\noindent{\bf Case 1:}  $f$ is a polynomial whose
coefficients have joint density. This is a consequence
of  Hammersley's formula, Theorem \ref{hammersley}.
%

\noindent{\bf Case 2:} $D$ is the unit disk or the
whole plane.

A Fubini argument
implies that there is a dense set $\mathcal R_D$ of
rectangles in $D$ such that for $R \in \mathcal R_D$,
almost surely $f$ does not vanish on $\partial R$.

It clearly suffices to show the claim when the
$\Lambda_i$ are disjoint elements of $\mathcal R_D$.

Let $\{f_M\}_{M \ge 1}$ denote the approximation of $f$ by
polynomials in Corollary \ref{approxcor}.

For $\Lambda \in \mathcal R_D$, the argument principle
implies that the number $N_M(\Lambda)$ of zeros of
$f_M$ in $\Lambda$, converges a.s.\ to the number
$N(\Lambda)$ of zeros of $f$ in $\Lambda$.

As $M$ varies, the random variables $\prod^n_{i=1}
N_M(\Lambda_i)$ are uniformly integrable, as they are
uniformly bounded in $L^2$ by Lemma \ref{combi}.

The covariance functions of $f_M$ and $f_M'$ converge
uniformly on each $\Lambda_i\times\Lambda_j$ to those
of $f$ and $f'$, whence the permanent-determinant formula (on
the right of \re{density}) for $f_M$ converges
uniformly on $\Lambda_1\times \ldots\times \Lambda_n$
to the permanent-determinant formula for $f$.

Applying formula \re{soshcor1} to $f_M$ and letting $M
\to \infty$ we see that it converges to the desired
formula for $f$.

\noindent{\bf Case 3:}  $D$ is simply connected: the
claim follows from the Riemann mapping theorem and
Case 2.

\noindent{\bf Case 4:} A general domain $D$. It
suffices to prove (\ref{soshcor1}) when
 each $\Lambda_j$ is a closed square in $D$.
Then we can find a simply connected domain $D_0
\subset D$ that contains all the $\Lambda_j$, and we
apply Case 3.
\end{proof}

For simple point processes, the following Lemma
implies Proposition \re{permformula} (ii).

\begin{lemma} \lb{combi}
Consider a simple point process (a random subset) $Z$
in a domain $D$ with counting function
$N(\Lambda)=\#(Z\cap \Lambda)$. Suppose that for any
disjoint Borel subsets $\Lambda_1,\dots ,\Lambda_k$ of
$D$, we have \bel{sosh} \ev \prod^k_{i=1} N(
\Lambda_i)= \int_{\Lambda_1\times\dots\times
\Lambda_k}
 p(z_1,\dots ,z_k)\,dz_1\dots dz_k \,.
\ee

Let $Z^{\wedge k}\subset Z^k$ denote the set of
$k$-tuples of distinct points. Then for any Borel set
$B \subset D^k $, we have
 \bel{soshcor2} \ev\,
\#(B\cap Z^{\wedge k}) =\int_Bp(z_1,\dots
,z_k)\,dz_1\dots dz_k \,. \ee
\end{lemma}
\begin{proof}
Note that both sides of \re{soshcor2} define a Borel
measure on subsets $B\in D^k$; thus it suffices to
show the equivalence for the case when $B=B_1\times
\cdots \times B_k$ is a product set.

Consider a finite Borel partition ${\mathcal P}$ of
$D$ and denote
$$
M_k({\mathcal P}) = \sum_{Q_1,\ldots,Q_k} \#(Q_1\times
\cdots \times Q_k\cap B \cap
Z^{k})=\sum_{Q_1,\ldots,Q_k} \prod_{i=1}^k N(Q_i \cap
B_i) \,
$$
where the sum is over ordered $k$-tuples
$(Q_1,\ldots,Q_k)$ of distinct elements of  ${\mathcal
P}$. Then the hypothesis (\ref{sosh}) implies that
\bel{ycollect} \ev M_k({\mathcal P})=
\sum_{Q_1,\ldots,Q_k}\int_{Q_1\times\dots\times Q_k
\cap B}
 p(z_1,\dots ,z_k)\,dz_1\dots dz_k\,,
\ee where we sum over the same $k$-tuples as above.
Consider a refining sequence of partitions ${\mathcal
P_j}$ of $D$ where the maximal diameter of the
elements of ${\mathcal P_j}$ tends to $0$ as $ j \to
\infty$. By definition, $M_k({\mathcal P})$ counts the
number of $k$-tuples $(z_1,\ldots,z_k)\in B$ where
$z_1,\ldots z_k$ are points of $Z$ in distinct
elements of $\mathcal P$. We deduce that
$$
M_k(\mathcal P_j) \to |B\cap Z^{\wedge k}|
$$
monotonically as $j \to \infty$. Taking expectations
and letting $j \to \infty$  yields (\ref{soshcor2}).
\end{proof}


We now proceed to analyze the behavior of Gaussian
analytic functions near their zeros.
\begin{lemma}\lb{simple}
Let $f$ be a Gaussian analytic function in a domain
$D$, and assume that for every $z\in D$ a.s.\
$z$ is not a double zero of $f$.
Then a.s.\ $f$ has no double zeros.
\end{lemma}
\noindent The Gaussian assumption is needed:
consider $(z-\gamma)^2$ with $\gamma$ a continuous random variable.
\begin{proof}
We may assume that there exists
$z_0 \in D$ such that $W=f(z_0)-\ev(f(z_0))$ is not identically zero
(otherwise there is nothing to show).
Let $g(z)=\ev(f(z)\overline{W})/\ev|W|^2$ and $h(z)=f(z)-Wg(z)$.
Then $g$ is a deterministic analytic function with $g(z_0)=1$,
 and $h(\cdot)$ is independent of $W$.
By assumption all the zeros of $g$ are not double zeros of $f$.
Any other double zero of $f$ would also be a double zero
of $\psi=W+h/g$. If $h/g$ is a random constant,
then $\psi$ a.s.\ has no zeros. Otherwise,
a.s.\ $\psi'=(h/g)'$ has at most countably many zeros
$\{\xi_j\}$, and they are a.s.\ not zeros of $\psi$
since
$W \ne -(h/g)(\xi_j)$ a.s.\ by independence.
\end{proof}

\begin{lemma}\lb{modulus}
Let $f$ be a Gaussian analytic function (not
necessarily mean 0) with radius of convergence $r_0$,
and let $M_r$ be its maximum modulus over the closed
disk of radius $r<r_0$. There exists $c,\gamma>0$ so
that for all $x>0$ we have
$$
\pr(M_r>x)\le c e^{-\gamma x^2}.
$$
\end{lemma}
\begin{proof}
Borell's Gaussian isoperimetric inequality (see
\cite{pollard}; the inequality was also shown
independently by \cite{tis75}) implies that for any
collection of mean 0 Gaussian variables with maximal
standard deviation $\sigma$, the max $M$ of the
collection satisfies \bel{borell}
\pr(M>\mbox{median}(M)+b\sigma) \le \pr(N>b) \, , \ee
where $N$ is standard normal. Now the median of $M_r$
is finite because $M_r<\infty$ a.s. Since the
distribution of $f(z)$ is continuous as a function of
$z$, the maximal standard deviation $\sigma$ in the
disk $B_r(0)$ is bounded. The mean 0 version of the
lemma follows by applying (\ref{borell}) to the real
and imaginary parts separately, and the general
version follows easily.
\end{proof}

\begin{lemma}\lb{diffzero}  Let $f(z)=a_0+a_1z+\ldots$ be a
Gaussian analytic function. Assume that $a_0$ is
nonconstant. Let $A_\eps$ denote the event that the
number of zeros of $f(z)$ in the disk $B_\eps$ about
$0$, differs from the number of zeros of
$h(z)=a_0+a_1z$ in $B_\eps$. \newline \noindent{\bf
(i)} For all $\delta>0$ there is $c>0$ (depending
continuously on the mean and covariance functions of
$f$) so that for all $\eps>0$ we have
$$
\pr(A_\eps)\le c\eps^{3-2\delta}.
$$
\noindent{\bf (ii)} $\pr(A_\eps \mid a_1, a_2,
\ldots)\le C \eps^{3}$, where $C$ may depend on
$(a_1,a_2,\ldots)$ but is finite almost surely.
\end{lemma}

\begin{proof} (i) By Rouch\'e's theorem, if $|h|>|f-h|$ on
the circle $\partial B_\eps$, then $f$ and $h$ have
the same number of zeros in $B_\eps$. By Lemma
\ref{modulus} applied to $(f-h)/z^2$,
 we have \bel{f-h}
 \pr(\max_{z\in\partial
B_\eps}|f(z)-h(z)|>\eps^{2-\delta})<c_0
\exp(-\gamma\eps^{-2\delta})< c_1 \eps^3. \ee
 Let $\Theta$ be the annulus $\partial
B_{a_1\eps}+B_{\eps^{2-\delta}}$, and consider the
following events:
\begin{eqnarray*}
D_0&=&\{|a_0|<2\eps^{1-\delta}\},\\
E&=&\{|a_1|<\eps^{-\delta}\},\\
F&=&\{\min_{z\in \partial B_\eps}
|h(z)|<\eps^{2-\delta}\} \, =\,\{-a_0\in \Theta\}.
\end{eqnarray*}
Note that $\pr E^c \le c_2\eps^{3}$ and that $E\cap
F\subset D_0$. Given $D_0$, the distribution of $a_0$
is approximately uniform on $B_{2\eps^{1-\delta}}$
(i.e., its conditional density is
$O(\eps^{2\delta-2})$).
 Since $\pr(E)$ tends to one as $ \eps \to 0$, this implies that
$$
\pr(F\cap E\gi D_0)=\pr(-a_0\in \Theta,E|D_0)\le c_3
\frac{\ev[\area(\Theta)]} {\area(B_{2\eps^{1-\delta}
}) } \le c_4 \eps^{2-\delta}/\eps^{1-\delta} =c_4\eps,
$$
and therefore
$$\pr(F)\le \pr(F\cap E|D_0)\pr(D_0)+\pr E^c\le
c_4\eps c_5 \eps^{2-2\delta}+c_2\eps^3\le
c_6\eps^{3-2\delta}.$$ Together with \re{f-h}, this
gives the desired result. Since all our bounds depend
continuously on the covariance function of $f$, we may
choose $c$ in a continuous manner, too.

\noindent{(ii)} The argument used to bound $\pr(F)$ in
$(i)$ also yields that
$$
\pr\Bigl(\min_{z\in \partial B_\eps}
|h(z)|<2|a_2|\eps^{2} \, \Big| \, \{a_j\}_{j \ge 1}
\Bigr) \le c_7 \eps^3 \,.
$$
An application of Rouch\'e's Theorem concludes the
proof.
\end{proof}


The following lemma relates the integral joint
intensity to the pointwise (strong) version.

\begin{lemma}\label{pwi}
Consider a simple point process on a domain $D$. Let
$z_j\in D$, $j=1,\ldots,n$. Assume that there exists
disjoint neighborhoods $D_j$ of $z_j$ and a $\delta>0$
so that the integral version of the joint intensity
satisfies
\begin{eqnarray*}
p(z_1,\ldots,z_n,z_*) <c_2 |z_j-z_*|^{-2+\delta}&&
\mbox{ on }D_1\times \cdots \times D_n\times
D_j\quad\mbox{ for all }j.
\end{eqnarray*}
Let $N_{j,\eps}$ denote the number of points in the
ball of radius $\eps$ about $z_j$. As $\eps\to 0$, we
have
 \bel{twointensities}\pr(N_{1,\eps}=\ldots
 =N_{n,\eps}=1)\le
\ev(N_{1,\eps}\cdots N_{n,\eps})=
\pr(N_{1,\eps}=\ldots =N_{n,\eps}=1)+o(\eps^{2n}).
 \ee
\end{lemma}

\begin{proof}
For nonnegative integers $N_j$ we have
 \bel{integers}
0\le \prod_{j=1}^n N_j - \prod_{j=1}^n \one(N_j = 1)
\le N_1\cdots N_n \sum_{k=1}^n (N_k-1).
 \ee
The left inequality is clear.  For the right one, if
for some $k$ we have $N_k>1$ then the $k$th term on
the right alone gives an upper bound. We apply
\re{integers} to the $N_{j,\eps}$ with small $\eps$
and take expectations. We apply Lemma \ref{combi} to
the set $B_\eps(z_1)\times \cdots \times
B_\eps(z_n)\times B_\eps(z_k)$ to get
$$
\ev(N_{1,\eps}\cdots N_{n,\eps}(N_{k,\eps}-1))
=\int_{B_\eps(z_1)\times \cdots \times
B_\eps(z_n)\times B_\eps(z_k)}p(w_1,\ldots, w_{n+1})
\,dw_1\cdots dw_{n+1}=o(\eps^{2n}).
$$
\end{proof}

\begin{lemma} \label{repell} Consider a Gaussian analytic function
$f$ in a domain $D$ with mean zero everywhere. Let
$z_1,\ldots, z_n\in D$, and assume that for each $j$,
the random variables $f'(z_j),$ $f(z_1),\ldots,$
$f(z_n)$ are linearly independent. Then there exists
neighborhoods $D_i$ of the $z_i$ so that for each
$1\le j \le n$, and $(w_1,\ldots, w_n,w_*)\in
D_1\times \cdots \times D_n\times D_j$ the integral
version of the joint intensity is defined and
satisfies
$$
p(w_1,\ldots, w_n,w_*) \le c |w_j-w_*|^2.
$$
\end{lemma}

\begin{proof}
By continuity, there exists bounded neighborhoods
$D_i$ of the $z_i$ so that

\begin{itemize}
\item[(i)] for all $w_*\in \bigcup D_i$ and
$w_i\in D_i$ the random variables $f'(w_*)$ and
$f(w_i)$, $1\le i\le n$ are linearly independent and
the determinant of the covariance matrix is bounded
away from $0$,

\item[(ii)] for all {\it distinct} points $w_*\in \bigcup D_i$ and
$w_i\in D_i$ the random variables $f(w_*)$ and
$f(w_i)$, $1\le i\le n$ are linearly independent.
\end{itemize}
Part (ii) follows by considering the Gaussian analytic
function $(f(w_*)-f(w_j))/(w_*-w_j)$, which has a removable singularity
at $w_j$.
Taylor expansion at $w$ implies that for  $w,z\in \bigcup D_i$, the
conditional distribution of $f'(w)$ given  $f(w)=f(z)=0$ is
Gaussian with variance bounded above by $c_1|z-w|^2$. Therefore,
$$
\ev \Bigl(|f'(w_1)\cdots f'(w_n)f'(w_*)|
\,\Big|\,f(w_1)=\ldots=f(w_*)=0 \Bigr) \le c_2 \eps^4,
$$
where $\eps$ is the distance between $w_*$ and the set
$\{w_1,\ldots ,w_n\}$. Furthermore,
\newcommand\dc{\operatorname{dc}}
$$\frac{\partial^2}{(\partial w_*)^2}
\det\Cov (f(w_1),\ldots,
f(w_n),f(w_*))\big\vert_{w_*=w_j} = \det
\Cov(f(w_1),\ldots , f(w_n),f'(w_j)).
$$
and since the right hand side is bounded away from
$0$, we get
$$|\det \Cov(f(w_1),\ldots f(w_n),f(w_*))| \ge c_3 \eps^2.$$

Now the permanent-determinant formula implies the
claim of the lemma.
\end{proof}

\begin{proofof}{Proposition \ref{permformula}} \

\noindent {\bf Step 1.} {We first verify the
equivalence of \re{prehan} and \re{density}.} First
note that when $f$ has Gaussian coefficients, then the
values and coefficients of $f$ are jointly Gaussian.
In particular, the expression \re{hamex} equals

%
\bel{recallmain} \begin{array}{l}\ev(|f'(z_1)\cdots
f'(z_n)|^2 \,\big|\,f(z_1)=\ldots=
f(z_n)=0)\,g(0,\ldots, 0) \\ \qquad = \ev\big(|\PP
f'(z_1)\cdots \PP f'(z_n)|^2)\,g(0,\ldots, 0),
 \end{array}\ee
  where $ g(0,\ldots,0) =
\pi^{-n} \det(A)^{-1} $ is the density of the Gaussian
vector $X=\big(f(z_j)\big)_{j=1}^n$ at 0, and $\PP$ is
the projection to the orthocomplement of the subspace
spanned by the entries of $X$.
Setting $Y=\big(f'(z_j)\big)_{j=1}^n$, note that the
projection $(I-\PP)$ of $Y$ onto the subspace spanned
by the entries of $X$ is given by $B A^{-1}X$. For a
column vector $Z$ of mean zero, recall that
$\Cov(Z)=\ev(ZZ^*)$. Now \bel{projcov} \Cov(\PP
Y)=\Cov(Y-  B A^{-1}X)=\Cov(Y)-\Cov(  B A^{-1}X)= C- B
A^{-1}  B^*. \ee
This proves the equivalence of \re{prehan} and
\re{density}.


\noindent{\bf Step 2.} Part (ii). By Lemma
\ref{simple}, the point process of zeros is simple.
Part (ii) follows from Proposition
\ref{prepermformula} and Lemma \ref{combi}.

\noindent{\bf Step 3.} Part (i). Let $F_\eps$ denote
the event that $f$ has a zero in each of
$B_\eps(z_1),\ldots,B_\eps(z_n)$. Since the function
$p(z_1,\ldots, z_n)$ is continuous, we have
\begin{eqnarray}\lb{pein}
\pr(F_\eps) &\le& \ev
N(B_\eps(z_1))\cdots N(B_\eps(z_n))
\\ \nonumber
&=&\int_{B_\eps(z_1)\times \cdots \times
B_\eps(z_n)}p(z_1,\ldots ,z_n)\,dz_1 \cdots dz_n
\\ \nonumber
&=& p(z_1,\ldots,z_n)\pi^n\eps^{2n}+o(\eps^{2n})
\end{eqnarray}
If, for some $j$, the derivative $f'(z_j)$ is not linearly
independent of $\{f(z_i) \, : \, 1\le i\le n\}$,
then $p(z_1,\ldots,z_n)=0$ and the claim follows from
\re{pein}. Otherwise, the claim follows from
\re{pein}, and Lemmas \ref{pwi} and \ref{repell}.
\end{proofof}

\smallskip

\noindent{\bf Acknowledgements.} We thank Russell
Lyons for useful discussions, Pavel Bleher and the referee for
references, and Alan Hammond and G\'abor Pete for
helpful remarks.
We are especially indebted to Ben Hough and Manjunath
Krishnapur for many insightful comments.
 Some of this work was written while the authors
participated in a workshop  at the Banff International Research
Station.
\bibliography{zeros}

\sc \bigskip \noindent Yuval Peres, Departments of
Statistics and Mathematics, U.C.\ Berkeley, CA 94720,
USA. \\{\tt peres@stat.berkeley.edu}, {\tt
stat-www.berkeley.edu/\~{}peres}

\sc \bigskip \noindent B\'alint Vir\'ag, Departments
of Mathematics and Statistics, University of Toronto,
ON, M5S 3G3, Canada. \\{\tt balint@math.toronto.edu},
{\tt www.math.toronto.edu/\~{}balint}


\end{document}